\documentclass[12pt]{article}
\usepackage{booktabs}
\usepackage{caption}
\usepackage{mathrsfs}
\usepackage{amsmath}
\usepackage{amsfonts,amsthm,amssymb,mathrsfs,bbding}
\usepackage{float}
\usepackage{graphicx}
\usepackage[hidelinks]{hyperref}
\usepackage{enumerate}
\usepackage{tikz}
\usepackage{caption}
\usepackage{rotating}
\allowdisplaybreaks[4]
\usepackage{tabularx}
\usepackage{cite}
\evensidemargin=\oddsidemargin\topmargin=-1.5cm

\newtheorem{thm}{Theorem}[section]
\newtheorem{prob}{Problem}[section]

\newtheorem{cor}{Corollary}[section]

\newtheorem{conj}{Conjecture}[section]

\theoremstyle{definition}

\addtocounter{section}{0}
\begin{document}

\title{{\bf Eigenvalues and factors: a survey} \footnote{This work was supported by the National Natural Science Foundation of China (Grant Nos. 12271162, 12301454, 12271425), Natural Science Foundation of Shanghai (Nos. 22ZR1416300 and 23JC1401500), the Program for Professor of Special Appointment (Eastern Scholar) at Shanghai Institutions of Higher Learning (No. TP2022031), and NRF-2021K2A9A2A1110161711.}}
\author{Dandan Fan$^{a,b}$,  Huiqiu Lin$^a$\thanks{Corresponding author;
Email addresses:
ddfan0526@163.com (D. Fan), huiqiulin@126.com (H. Lin), luhongliang@mail.xjtu.edu.cn (H. Lu), suil.o@sunykorea.ac.kr (O. S).}, Hongliang Lu$^c$,  Suil O$^d$\\[2mm]
\small\it $^a$ School of Mathematics, East China University of Science and Technology, \\
\small\it   Shanghai 200237, China\\[1mm]
\small\it $^b$ College of Mathematics and Physics, Xinjiang Agricultural University\\
\small\it Urumqi, Xinjiang 830052, China\\[1mm]
\small\it $^c$ School of Mathematics and Statistics, Xi'an Jiaotong University\\
\small\it Xi'an, Shaanxi 710049, China\\[1mm]
\small\it $^d$ Department of Applied Mathematics and Statistics, The State University of New York\\
\small\it Korea, Incheon, 21985, Republic of Korea
}
\date{}
\maketitle
{\flushleft\large\bf Abstract}
A factor of a graph is a spanning subgraph satisfying some given conditions. An earlier survey of factors can be traced back to the Akiyama and Kano [J. Graph Theory, 1985, 9: 1-42] in which they described the characterization of factors in (bipartite) graphs and digraphs, respectively. Soon after, Kouider and Vestergaard summarized the findings related to  connected factors [Graphs Combin., 2005, 21(1): 1-26]. Plummer extended the aforementioned research by providing a comprehensive overview of progress made in the study of graph factors and factorization from 1985 to 2003 [Discrete Math., 2007, 7-8: 791-821]. In this paper, we aim to summarize the relevant results regarding factors from the perspective of eigenvalues.

\begin{flushleft}
\textbf{Keywords:} eigenvalue; factor; perfect matching; Hamiltonian cycle.
\end{flushleft}
\textbf{AMS Classification:} 05C50

\newpage

\tableofcontents

\section{Introduction}

The initial study of factors was due to Danish mathematician Petersen in 1891, who proved that every graph of even degrees can be decomposed into the union of edge-disjoint 2 factors. Additionally, Petersen also demonstrated that every 2-connected cubic graph possesses a 1-factor. These two results can be regarded as precursors to factor theory.

Based on the properties of induced subgraphs, the factor problem can be divided into two classes: degree-constrained factors and component factors. Degree-constrained factors refer to factors where all vertices have degrees that do not exceed a given value. Common types of degree-constrained factors include odd factors, even factors, $k$-factors and $[a,b]$-factors, and so on. These two factors are not completely disjoint. For instance, a 1-factor in a graph can be regard as a spanning subgraph where each vertex has the same degree 1, or as a spanning subgraph where each component is a path $P_2$. Factor problems may also overlap with other graph theory problems, such as Hamiltonian cycle problem and subgraph packing problem.

For a graph $G$, let $A(G)$ denote  the adjacency matrix of $G$ and let $\lambda_i(G)$ denote the $i$th largest eigenvalue of $A(G)$. Particularly, the largest eigenvalues of $A(G)$, denoted by $\rho(G)$, is called the \textit{spectral radius} of $G$. The \textit{Laplacian matrix} of $G$ is defined as $L(G)=D(G)-A(G)$, where $D(G)$ is the diagonal matrix of vertex degrees of $G$. The Laplacian matrix is positive semidefinite and we order its eigenvalues as $\mu_{1}\geq\mu_{2}\geq\ldots\geq\mu_{n-1}\geq\mu_{n}=0$.

In this paper, we only consider simply graphs unless otherwise stated. The paper is organized as follows. In section 2, we will provide an overview of eigenvalue conditions, specifically highlighting the spectral radius and third largest eigenvalue, which guarantee the existence of a degree-constrained factor in graphs. In section 3, we summarize relevant eigenvalue results concerning the existence of component factors such as path factors, star factors, Hamiltonian (path) cycles, and $k$-trees. In section 4, we conclude the eigenvalue and factor packing problem, with a specific focus on 1-factors and edge-disjoint spanning trees. Moreover, we list some problems in the end.


\section{Eigenvalues and degree-constrained factors}
In this section, we focus on the eigenvalue conditions for the existence of degree-constrainted factors, which include many well-known factors, such as 1-factors, $k$-factors and $[a,b]$-factors, and so on.
\subsection{Matchings}
\subsubsection{Matchings and perfect matchings}
A \textit{matching} in a graph is a set of pairwise nonadjacent edges. The matching number $\beta(G)$ of a graph $G$ is the cardinality of a maximum matching of $G$. Erd\H{o}s and Gallai\cite{Erdos-Gallai} determined the maximum number of edges in a graph of order $n$ with given matching number. Let $G_{1} \nabla G_{2}$ be the graph obtained from the disjoint union $G_{1}\cup G_{2}$ by adding all edges between $G_{1}$ and $G_{2}$. Denote by $e(G)$ the number of edges in $G$.
\begin{thm}[Erd\H{o}s and Gallai\cite{Erdos-Gallai}]
Suppose that $G$ is a graph of order $n$ with matching number $\beta$. If $n > 2\beta$, then
$$e(G) \leq \max \Big\{{2\beta+1\choose 2}, {\beta\choose 2}+\beta(n-\beta)\Big\}.$$
If $2 \beta+1 \leq n<(5\beta+3)/2$, then $K_{2\beta+1} \cup (n-2\beta-1)K_1$ is the unique extremal graph; if $n=(5\beta+3)/2$, then there are two extremal graphs $K_{2\beta+1} \cup (n-2\beta-1)K_1$ and $K_{\beta} \nabla (n-\beta)K_1$. If $n>(5\beta+3)/2$, then $K_{\beta} \nabla (n-\beta)K_1$ is the unique extremal graph.
\end{thm}

Let $o(G)$ be the number of odd components in a graph $G$. Berge \cite{Berge59} gave the formula, which is called by Berge-Tutte Formula.
\begin{thm}[Berge \cite{Berge59}]
For a graph $G$ of order $n$,
\[
\beta(G)=\frac{1}{2}(n-\max_{S\subseteq V(G)}\{o(G-S)-|S|\}).
\]
\end{thm}

Brouwer and  Gregory\cite{SG07} utilized the $(k+1)$th largest adjacency eigenvalue to provide a lower bound on the matching number in a regular graph. For a positive integer $r$, a graph is $r$-regular if every vertex has the degree $r$.
\begin{thm}[Brouwer and  Gregory\cite{SG07}]
Let $G$ be a connected $r$-regular graph of even order $n$. Suppose that $k>0$ is an integer such that $n\equiv k\pmod 2$. If
$$
r-\lambda_{k+1}(G)> \left\{
\begin{array}{ll}
 0.1457& \mbox{if $r=3$,}\vspace{1ex}\\
1-\frac{3}{r+1}- \frac{1}{(r+1)(r+2)}& \mbox{if $r$ is even,}\\
1-\frac{4}{r+2}- \frac{1}{(r+2)^2}& \mbox{if $r\geq 5$ is odd.}
\end{array}
\right.
$$
then $\beta(G)\geq (n-k)/2+1$.
\end{thm}

Based on Berge-Tutte Formula, Feng, Yu and Zhang\cite{Feng-Yu-Zhang} provided a spectral radius condition with given matching number.
\begin{thm}[Feng, Yu and Zhang\cite{Feng-Yu-Zhang}]
Let $G$ be a graph of order $n$ with matching number $\beta$.
\begin{enumerate}[(i)]
\item If $n=2 \beta$ or $2 \beta+1$, then $\rho(G) \leq \rho\left(K_n\right)$, with equality if and only if $G\cong K_n$.
\item If $2 \beta+2 \leq n<3 \beta+2$, then $\rho(G) \leq 2 \beta$, with equality if and only if $G\cong K_{2 \beta+1} \cup (n-2 \beta-1)K_1$.
\item If $n=3 \beta+2$, then $\rho(G) \leq 2 \beta$, with equality if and only if $G\cong K_\beta \nabla (n-\beta)K_1$ or $G\cong K_{2 \beta+1} \cup (n-2 \beta-1)K_1$.
\item If $n>3 \beta+2$, then $\rho(G) \leq (\beta-1+\sqrt{(\beta-1)^2+4 \beta(n-\beta)})/2$, with equality if and only if $G\cong K_\beta \nabla (n-\beta)K_1$.
\end{enumerate}
\end{thm}
Recognizing that the graphs considered by Feng, Yu, and Zhang may not necessarily be connected, Chen and Lu\cite{Chen-Lu} undertaken the following research efforts specifically focusing on connected graphs.
\begin{thm}[Chen and Lu\cite{Chen-Lu}]
Let $G$ be a connected graph of order $n$ with matching number $\beta$.
\begin{enumerate}[(i)]
\item  If $n=2\beta$ or $2\beta+1$, then $\rho(G) \leq \rho(K_n)$, with equality if and only if $G\cong K_n$.
\item  If $2\beta+2 \leq n \leq 3\beta-1$, then $\rho(G) \leq \rho(K_1 \nabla(K_{2\beta-1} \cup(n-2\beta) K_1))$, with equality if and only if $G\cong K_1 \nabla(K_{2\beta-1} \cup(n-2\beta) K_1)$.
\item  If $n \geq 3\beta$, then $\rho(G)\leq (\beta-1+\sqrt{(\beta-1)^2+4 \beta(n-\beta)})/2$, with equality if and only if $G \cong K_{\beta}\nabla (n-\beta)K_1$.
\end{enumerate}
\end{thm}


 Kim, O, Sim and Shin\cite{Kim-O-Sim-Shin} proved an upper bound for the spectral radius conditions in an $n$-vertex connected graph $G$ with $\beta(G) \leq\frac{n-k}{2}$.
\begin{thm}[Kim, O, Sim and Shin\cite{Kim-O-Sim-Shin}]
Let $n$ and $k$ be two positive integers, and let $G$ be a connected graph of order $n$ with matching number $\beta \leq\frac{n-k}{2}$, where $2 \leq k \leq n-2$ and $n \equiv k(\bmod~2)$. Then we have the following conclusions.
\begin{enumerate}[(i)]
\item If $n \leq 3 k$, then
$$\rho(G) \leq \rho\Big(K_{\frac{n-k}{2}} \nabla \frac{n+k}{2}K_1\Big),$$
 with equality if and only if $G\cong K_{\frac{n-k}{2}} \nabla \frac{n+k}{2}K_1$.
\item If $n \geq 3 k+2$, then
$$\rho(G) \leq \rho(K_1 \nabla(K_{n-k-1} \cup kK_1)),$$
 with equality if and only if $G\cong K_1 \nabla(K_{n-k-1} \cup kK_1)$.
\end{enumerate}
\end{thm}

For $\beta<\lfloor\frac{n}{2}\rfloor$, Zhang, Wang and Wang\cite{Zhang-Wang-Wang} utilized Berge-Tutte Formula to obtain that $\delta(G) \leq \beta$. Take the minimum degree into consideration, they characterized the extremal graphs having maximum spectral radius with fixed matching number.
For a graph $H$, let $K_1 \nabla_\delta K_a \nabla H$ be the graph obtained from $K_a \nabla H$ and a new vertex $w$ by connecting $w$ to (any) $\delta$ vertices of the part $K_a$ in $K_a \nabla H$, where $1 \leq \delta \leq a$. Denote by $\delta(G)$ the minimum degree of a graph $G$. Clearly, $\delta(K_1 \nabla_\delta K_a \nabla H)=\delta$.
\begin{thm}[Zhang, Wang and Wang\cite{Zhang-Wang-Wang}]
Suppose that $n, k, \delta$ are three positive integers, where $n \geq 29, 2\leq k \leq n-2, n \equiv$ $k(\bmod~2)$, and $1 \leq \delta \leq \frac{n-k}{2}$. Let $G$ be a connected graph of order $n$ with minimum degree $\delta$ and matching number $\beta \leq \frac{n-k}{2}$. Then
$$
\rho(G) \leq \max\Big\{\rho(K_\delta \nabla(K_{n+1-2 \delta-k} \cup (\delta+k-1)K_1), \rho\Big(K_1 \nabla_\delta K_{\frac{n-k}{2}} \nabla \Big(\frac{n+k}{2}-1\Big)K_1\Big)\Big\},
$$
with equality if and only if $G\cong K_\delta \nabla(K_{n+1-2 \delta-k} \cup (\delta+k-1)K_1$ or $G\cong K_1 \nabla_\delta K_{\frac{n-k}{2}} \nabla $ $(\frac{n+k}{2}-1)K_1$.
\end{thm}
Let
$$
\Delta(n, t, k)=(n^2-k^2)^2-2(n^2-k^2)(k+t-1)(n+k+10 t-4)+16 t(k+t-1)^2(n+4t-k-2),
$$
where $n, k$ and $t$ are three positive integers.
\begin{thm}[Zhang\cite{Zhang2022}]
Suppose that $n, t$ and $k$ are three positive integers, where $2 \leq k \leq n-2,1 \leq t \leq \frac{n-k}{2}$ and $n \equiv k(\bmod~2)$. Let $G$ be a t-connected graph of order $n$ with matching number $\beta\leq \frac{n-k}{2}$. Then we have the following conclusions.
\begin{enumerate}[(i)]
\item If $\Delta(n, t, k)>0$, then
$$\rho(G) \leq \rho(K_t \nabla (K_{n+1-2 t-k} \cup (t+k-1)K_1)),$$
with equality if and only if $G\cong K_t \nabla (K_{n+1-2 t-k} \cup (t+k-1)K_1)$.
\item If $\Delta(n, t, k)=0$, then
$$\rho(G) \leq \rho(K_t \nabla (K_{n+1-2 t-k} \cup (t+k-1)K_1))=\rho\Big(K_{\frac{n-k}{2}} \nabla \frac{n+k}{2}K_1\Big),$$ with equality if and only if $G\cong K_t \nabla(K_{n+1-2 t-k} \cup (t+k-1)K_1)$ or $G\cong K_{\frac{n-k}{2}} \nabla \frac{n+k}{2}K_1$.
\item If $\Delta(n, t, k)<0$, then
$$\rho(G) \leq \rho\Big(K_{\frac{n-k}{2}} \nabla \frac{n+k}{2}K_1\Big),$$
 with equality if and only if $G\cong K_{\frac{n-k}{2}} \nabla \frac{n+k}{2}K_1$.
\end{enumerate}
\end{thm}

Let $\mathcal{F}$ be a family of graphs. A graph $G$ is called \textit{$\mathcal{F}$-free} if it does not contain any graph in $\mathcal{F}$ as a subgraph. The classical Brualdi-Solheid-Tur\'{a}n type problems consider the maximum spectral radius of an $n$-vertex $\mathcal{F}$-free graph. Let $T_k(n)$ be the complete $k$-partite graph of order $n$ whose partition sets have sizes as equal as possible. In 2007, Nikiforov\cite{Nikiforov2007} characterized the extremal graph among all $K_{r+1}$-free graphs.
\begin{thm}[Nikiforov\cite{Nikiforov2007}]
If $G$ is a $K_{r+1}$-free graph of order $n$, then $\rho(G) \leq \rho(T_k(n))$, with equality if and only if $G \cong T_k(n)$.
\end{thm}
Let $\operatorname{ex}(n,\mathcal{F})$ and $\operatorname{spex}(n,\mathcal{F})$ be the maximum number of edges and maximum value of spectral radius among all $\mathcal{F}$-free graphs of order $n$, respectively. In the meantime, Feng, Yu and Zhang\cite{Feng-Yu-Zhang} determined the exact value of $\operatorname{spex}\left(n, M_{s+1}\right)$, where $M_{s+1}$ denotes a matching with $s+1$ edges. Alon and Frankl\cite{Alon-Frankl} provided further research by showing that $\operatorname{ex}(n,\{K_{k+1}, M_{s+1}\})=\max \{e(T_k(2 s+1)),e(G_k(n, s))\}$, where $G_k(n, s)=T_{k-1}(s) \nabla (n-s)K_1$. Later, Wang, Hou and Ma\cite{Wang-Hou-Ma} determined the exact value of $\operatorname{spex}(n,\{K_{k+1}, M_{s+1}\})$.
\begin{thm}[Wang, Hou and Ma\cite{Wang-Hou-Ma}]
For $n \geq 4 s^2+9 s$ and $k \geq 2$,
$$
\operatorname{spex}(n,\{K_{k+1}, M_{s+1}\})=\rho(G_k(n, s)).
$$
\end{thm}

A \textit{perfect matching} is a matching that covers all vertices of the graph. It is well known that $\beta(G) = \frac{n}{2}$ if $G$ has a perfect matching. Perfect matching as the simplest degree-constrained factor in a graph has been widely investigated. Obviously, a perfect matching is a 1-factor of graphs. Tutte\cite{Tutte1947} gave the characterization of 1-factors in 1947, which remains one of the cornerstone results in factor theory. Recall that $o(G)$ is the number of odd components in a graph $G$.
\begin{thm}[Tutte\cite{Tutte1947}]
A graph $G$ has a perfect matching if and only if
$$o(G-S)\leq |S|$$
for all $S\subseteq V(G)$.
\end{thm}

Utilizing this Theorem, Brouwer and Haemers\cite{A.B} described a regular graph to contain a perfect matching in terms of the third largest eigenvalue in 2005.

\begin{thm}[Brouwer and Haemers\cite{A.B}]\label{Brouwer and Haemers}
Let $G$ be a connected $r$-regular graph of even order $n$. If
$$
\lambda_{3}(G)\leq \left\{
\begin{array}{ll}
r-1+\frac{3}{r+1} & \mbox{if $r$ is even,}\vspace{1ex}\\
r-1+\frac{3}{r+2} & \mbox{if $r$ is odd,}
\end{array}
\right.
$$
then $G$ has a perfect matching.
\end{thm}
In \cite{Cioaba2005}, for $r$ is odd, the upper bound in Theorem \ref{Brouwer and Haemers} can be improved to $\lambda_3 \leq r- 1+\frac{4}{r+2}$. Later, Cioab\u{a}, Gregory and Haemers\cite{S.C} further improved the above upper bounds as follows.
\begin{thm}[Cioab\u{a}, Gregory and Haemers\cite{S.C}]\label{thm::1.2}
Let $G$ be a connected $r$-regular graph of even order $n$. If
$$
\lambda_{3}(G)< \left\{
\begin{array}{ll}
\theta=2.85577\ldots & \mbox{if $r=3$,}\vspace{1ex}\\
\frac{r-2+\sqrt{r^{2}+12}}{2} & \mbox{if $r\geq 4$ is even,}\vspace{1ex}\\
\frac{r-3+\sqrt{(r+1)^{2}+16}}{2} & \mbox{if $r\geq 5$ is odd,}\\
\end{array}
\right.
$$
where $\theta$ is the largest root of $x^{3}-x^{2}-6x+2=0$, then $G$ has a perfect matching.
\end{thm}

The graphic parameters play an important role in the study of factors. By imposing the minimum degree of a graph as a parameter, Liu, Liu and Feng\cite{Liu-Liu-Feng} extended the above results to graphs that are not necessarily regular.
\begin{thm}[Liu, Liu and Feng\cite{Liu-Liu-Feng}]\label{Liu-Liu-Feng}
Let $G$ be a connected graph of even order $n$ with minimum degree $\delta\geq 2$. If $n\geq \max\{7+7\delta+2\delta^{2}, \delta^{3}+3\delta^{2}+2\delta\}$ and
$$\rho(G)\geq\rho(K_{\delta}\nabla(K_{n-2\delta-1}\cup (\delta+1)K_1)),$$
then $G$ has a perfect matching unless $G\cong K_{\delta}\nabla (K_{n
-2\delta-1}\cup(\delta+1)K_1)$.
\end{thm}

In 2021, O\cite{S.O} provided an edge condition to guarantee the existence of a perfect matching in a connected graph for general $n$.
\begin{thm}[O\cite{S.O}]\label{edge}
Let $G$ be a connected graph of even order $n$.
\begin{enumerate}[(i)]
\item For $n \geq 10$ or $n=4$, if $e(G)>{n-2\choose 2}+2$, then $G$ has a perfect matching.
 \item For $n=6$ or $n=8$, if $e(G)>9$ or $e(G)>18$, then $G$ has a perfect matching.
 \end{enumerate}
\end{thm}

Observe that none of $K_{1}\nabla (K_{n-3}\cup 2K_{1})$, $K_2\nabla 4K_1$ and $K_3\nabla 5K_1$ contains a perfect matching. Combining this with $e(K_{1}\nabla (K_{n-3}\cup 2K_{1}))={n-2\choose 2}+2$ for $n\geq 4$, $e(K_2\nabla 4K_1)=9$ and $e(K_3\nabla 5K_1)=18$, it is easy to find that the bounds in Theorem \ref{edge} are sharp. Moreover, O\cite{S.O} obtained a similar result from the perspective of spectral radius.

\begin{thm}[O\cite{S.O}]\label{spectral}
Let $G$ be a connected graph of even order $n$.
\begin{enumerate}[(i)]
\item For $n\geq 8$ or $n=4$, if $\rho(G)>\theta(n)$ where $\theta(n)$ is the largest root of $x^3-(n-4)x^2-(n-1)x+2(n-4)=0$, then $G$ has a perfect matching.
\item For $n=6$, if $\rho(G)>\frac{1+\sqrt{33}}{2}$, then $G$ has a perfect matching.
\end{enumerate}
\end{thm}

 One can verify that $\rho(K_{2}\nabla 4K_{1})=\frac{1+\sqrt{33}}{2}$ and $\rho(K_{1}\nabla (K_{n-3}\cup 2K_{1}))=\theta(n)$. This implies that the conditions in Theorem \ref{spectral} are best possible.

For matchings in bipartite graphs, the K\"{o}nig-Hall Theorem, also known as Hall's Theorem or the Marriage Theorem, was established by K\"{o}nig (1931) and Hall (1935). Due to its wide applications to many graph theory problems and to other branches of mathematics, K\"{o}nig-Hall Theorem remains one of the most influential graph-theoretic results. For any $v\in V(G)$, let $N_{G}(v)$ denote the neighborhood of $v$ in $G$, and for $S\subseteq V(G)$, let $N_{G}(S)=\bigcup_{v\in S}N_{G}(v)$.

\begin{thm}[Hall \cite{Hall}]\label{Hall}
A bipartite graph $G=(X,Y)$ has a perfect matching if and only if $|X|=|Y|$ and
$$|N_{G}(S)|\geq|S|$$
for all $S\subseteq X$.
\end{thm}

A bipartite graph $G = (X, Y)$ is called \emph{balanced} if $|X| = |Y|$. Motivated by Hall's condition, Fan, Goryainov, Huang and Lin\cite{D.F} gave a spectral radius condition for the existence of a perfect matching in a balanced bipartite graph  with fixed minimum degree.
Given two bipartite graphs $G_1=(X_{1}, Y_{1})$ and $G_2=(X_{2}, Y_{2})$, let $G_{1}\nabla_{1} G_2$ denote the graph obtained from $G_1 \cup G_2$ by adding all possible edges between $X_{2}$ and $Y_{1}$.

\begin{thm}[Fan, Goryainov, Huang and Lin\cite{D.F}]\label{Fan-Goryainov-Huang-Lin}
Let $G$ be a balanced bipartite graph of order $2n$  with minimum degree $\delta\geq 1$. If
 $$\rho(G)\geq\rho(K_{\delta+1,\delta}\nabla_{1} K_{n-\delta-1,n-\delta}),$$
then $G$ has a perfect matching unless $G\cong K_{\delta+1,\delta}\nabla_{1} K_{n-\delta-1,n-\delta}$.
\end{thm}

By a simple computation, it is easy to find that
$$\rho(K_{\delta+1,\delta}\nabla_{1} K_{n-\delta-1,n-\delta})=\frac{\sqrt{2n^2 - 2(\delta + 1)n +  2\delta(\delta+1) + 2f(n,\delta)}}{2},$$
 where $f(n,\delta)=\sqrt{n^4 - 2(\delta + 1)n^3 + (1 - \delta^2)n^2 + 2\delta(3\delta^2 + 4\delta + 1)n - 3\delta^2(\delta+1)^2}$. Considering that the expression of $\rho(K_{\delta+1,\delta}\nabla_{1} K_{n-\delta-1,n-\delta})$  is complicated and the fact that $\rho(K_{\delta+1,\delta}\nabla_{1} K_{n-\delta-1,n-\delta})\geq \rho(K_{n-\delta-1,n}\cup (\delta+1)K_1)=\sqrt{n(n-\delta-1)}$,
  they provided the following result which improves Theorem \ref{Fan-Goryainov-Huang-Lin} for sufficiently large $n$.

\begin{thm}[Fan, Goryainov, Huang and Lin\cite{D.F}]
Let $G$ be a balanced bipartite graph of order $2 n$ with minimum degree $\delta \geq 1$. If $n \geq \frac{1}{2} \delta^3+\frac{1}{2} \delta^2+\delta+4$ and $$\rho(G) \geq \sqrt{n(n-\delta-1)},$$
then $G$ has a perfect matching unless $G \cong K_{\delta+1, \delta} \nabla_1 K_{n-\delta-1, n-\delta}$.
\end{thm}

\subsubsection{Rainbow matchings}
Let $\mathcal{G}=\{G_1,\ldots,G_t\}$ be a family of not necessarily distinct graphs with common vertex set $V$. We say that a graph $H$ with vertex set $V$ is \emph{$\mathcal{G}$-rainbow} if there exists a bijection $\phi: E(H)\rightarrow[t]$ such that $e\in E(G_{\phi(e)})$ for each $e\in E(H)$. Guo, Lu, Ma and Ma \cite{GLMM} generalized the results in \cite{Feng-Yu-Zhang, Chen-Lu} to rainbow version.

\begin{thm}[Guo, Lu, Ma and Ma\cite{GLMM}]
Suppose that $n$ and $m$ are two positive integers such that $n\geq 2m+2$. Let $\mathcal{G}=\{G_1,G_2,\ldots,G_{m+1}\}$ be a family of  graphs on the same vertex set $[n]$. If
\begin{align*}
 \rho(G_i)\geq\max\Big\{2m,\frac{1}{2}(m-1+\sqrt{(m-1)^2+4m(n-m)})\Big\},
\end{align*}
 then
	\begin{enumerate}[(i)]
		\item  for $2m+2\leq n<3m+2$, $\mathcal{G}$ admits a rainbow matching unless $G_1=\ldots=G_{m+1}\cong K_{2m+1}\cup (n-2m-1)K_1$;
		\item  for $n=3m+2$, $\mathcal{G}$ admits a rainbow matching unless $G_1=\ldots=G_{m+1}\cong K_{m}\nabla(n-m)K_1$ or $G_1=\ldots=G_{m+1}\cong K_{2m+1}\cup (n-2m-1)K_1$;
		\item  for $n>3m+2$, $\mathcal{G}$ admits a rainbow matching unless $G_1=\ldots=G_{m+1}\cong K_{m}\nabla(n-m)K_1$.
	\end{enumerate}
\end{thm}

\subsubsection{Fractional matchings and fractional perfect matchings}
A \textit{fractional matching} of a graph $G$ is a function $f$ giving each edge a number in $[0,1]$ such that $\sum_{e\in \Gamma(v)}f(e)\leq 1$ for each $v\in V(G)$, where $\Gamma(v)$ is the set of edges incident to $v$. The fractional matching number of $G$, denote by $\alpha_*^{\prime}(G)$, is the maximum of $\sum_{e \in E(G)} f(e)$ over all fractional matchings $f$.

Let $i(G)$ denote the number of isolated vertices of $G$. By replacing $o(G-S)$ with $i(G-S)$ in graph $G-S$, Scheinerman and Ullman\cite{Scheinerman-Ullman} provided a fractional version of the Berge-Tutte Formula, that is, $\alpha_*^{\prime}(G)=\min_{S\subseteq V(G)}\frac{n-i(G-S)+|S|}{2}$.
In \cite{O2016}, O established an upper bound for the spectral radius of a graph $G$ of order $n$ with minimum degree $\delta$, stating that $\rho(G)<\delta\sqrt{1+2k/(n-k)}$, which guarantees that its fractional matching number satisfies $\alpha_*^{\prime}(G)>\frac{n-k}{2}$ where $k$ is a real number
between 0 and $n$. In the same paper, O further investigated the relationships between the $\rho(G)$, minimum degree $\delta$ and $\alpha_*^{\prime}(G)$.
For two positive integers $\delta$ and $k$, let $\mathcal{H}(\delta, k)$ be the family of connected bipartite graphs $H$ with the bipartitions $A$ and $B$ such that $|A|=|B|+k$, $d_{H}(v)=\delta$ for $v\in A$ and the degrees of vertices in $B$ are equal.
\begin{thm}[O\cite{O2016}]
If $G$ is a graph of order $n$ with minimum degree $\delta$, then
$$
\alpha_*^{\prime}(G) \geq \frac{n \delta^2}{\rho(G)^2+\delta^2},
$$
with equality if and only if $k=\frac{n(\rho(G)^2-\delta^2)}{\rho(G)^2+\delta^2}$ is an integer and $G$ is an element of $\mathcal{H}(\delta, k)$.
\end{thm}

        Luo, Liu and Ao\cite{Luo-Liu-Ao} considered a lower bound of spectral radius in a graph $G$ to guarantee $\alpha_*^{\prime}(G)>\frac{n-k}{2}$.
\begin{thm}[Luo, Liu and Ao\cite{Luo-Liu-Ao}]
Let $k$ and $n$ be two positive integer with $k<n$, and let $G$ be a graph of order $n \geq \max \left\{6 \delta+5 k+1,5 \delta+k^2+4 k+1\right\}$ with minimum degree $\delta$. If
$$
\rho(G) \geq \rho(K_\delta \nabla(K_{n-2 \delta-k}\cup (\delta+k) K_1)),
$$
then $\alpha_*^{\prime}(G)>\frac{n-k}{2}$ unless $G \cong K_\delta \nabla(K_{n-2 \delta-k}\cup (\delta+k) K_1)$.
\end{thm}

A  \textit{fractional perfect matching} of an $n$-vertex graph $G$ is a fractional matching $f$ with $\sum_{e\in E(G)}f(e)=\frac{n}{2}$.
\begin{thm}[Scheinerman and Ullman\cite{Scheinerman-Ullman}]\label{Scheinerman-Ullman}
A graph $G$ has a fractional perfect matching if and only if
$$i(G-S)\leq |S|$$
for all subset $S\subseteq V(G)$.
\end{thm}

Given a graph $G$, if it does not have fractional perfect matchings, by Theorem \ref{Scheinerman-Ullman}, there exists a vertex set $S_0\subseteq V(G)$ such that $i(G-S_0)>|S_0|$. Along this line, Fan, Lin and Lu\cite{Fan-Lin-Lu} proved the following result.
\begin{thm}[Fan, Lin and Lu\cite{Fan-Lin-Lu}]
Suppose that $G$ is a connected graph of order $n$ with minimum degree $\delta$. If $n\geq 8\delta+4$ and $$\rho(G)\geq\rho(K_{\delta} \nabla (K_{n-2\delta-1}\cup (\delta+1)K_{1})),$$
 then $G$ has a fractional perfect matching unless $G\cong K_{\delta} \nabla (K_{n-2\delta-1}\cup (\delta+1)K_{1})$.
\end{thm}

In \cite{Pan-Liu,Li-Miao-Zhang}, the spectral radius conditions of the existence of perfect fractional matching in a graph for general $n$ were also provided.

%

\subsection{$[1,b]$-factors}

An $\textit{odd [1,b]-factor}$ of a graph $G$ is a spanning subgraph $H$ such that $d_{H}(v)$ is odd and $1\leq d_{H}(v)\leq b$ for each $v\in V(G)$. As an extension of Tutte's 1-Factor Theorem, the well-known sufficient and necessary condition for the existence of an odd $[1,b]$-factor was established by Amahashi.
\begin{thm}[Amahashi\cite{A.A}]\label{Amahashi}
Let $G$ be a graph and let $b$ be a positive odd integer. Then $G$ has an odd $[1,b]$-factor if and only if
$$o(G-S)\leq b|S|$$
for all $S\subseteq V(G)$.
\end{thm}
For $X,Y\subseteq V(G)$, we denote by $E_{G}(X,Y)$ the set of edges with one endpoint in $X$ and one endpoint in $Y$, and $e_{G}(X,Y)=|E_{G}(X,Y)|$.
Theorem \ref{Amahashi} states that there exists a subset $S \subseteq V(G)$ such that $o(G-S)>b|S|$ if there is no odd $[1, b]$-factor in an $r$-regular graph $G$. By counting the number of edges between $S$ and $G-S$, Lu, Wu and Yang \cite{Lu-Wu-Yang} showed that $G-S$ has at least three odd components $G_1, G_2, G_3$ such that $e_{G}(V(G_i), S)<\lceil\frac{r}{b}\rceil$ for $1\leq i\leq 3$. Then they gave a sufficient condition for the existence of an odd $[1,b]$-factor in a graph in terms of the third largest eigenvalue.
\begin{thm}[Lu, Wu and Yang\cite{Lu-Wu-Yang}]\label{Lu-Wu-Yang}
Let $G$ be a connected $r$-regular graph of even order, where $r\geq 3$. If
$$
\lambda_{3}(G)\leq \left\{
\begin{array}{ll}
 r-\frac{\lceil\frac{r}{b}\rceil-2}{r+1}+\frac{1}{(r+1)(r+2)} & \mbox{ if $r$ is even and $\lceil\frac{r}{b}\rceil$ is even,}\\
 r-\frac{\lceil\frac{r}{b}\rceil-1}{r+1}+\frac{1}{(r+1)(r+2)} & \mbox{if $r$ is even and $\lceil\frac{r}{b}\rceil$ is odd,}\\
r-\frac{\lceil\frac{r}{b}\rceil-1}{r+2}+\frac{1}{(r+2)^{2}} & \mbox{if $r$ is odd and $\lceil\frac{r}{b}\rceil$ is even,}\\
 r-\frac{\lceil\frac{r}{b}\rceil-2}{r+2}+\frac{1}{(r+2)^{2}} & \mbox{if $r$ is odd and $\lceil\frac{r}{b}\rceil$ is odd},
\end{array}
\right.
$$
then $G$ has an odd $[1,b]$-factor.
\end{thm}

They found the upper bounds appeared above are reachable in the family $\mathcal{F}_{r, b}$, where $\mathcal{F}_{r, b}$ is a family of such a possible component depending on $r$ and $b$. Later, Kim, O, Park and Ree\cite{Kim-O-Park-Ree} improved the spectral conditions of Theorem \ref{Lu-Wu-Yang} by constructing a graph $H_{r,\eta}$ in $\mathcal{F}_{r, b}$ with $\lambda_1(H_{r,\eta})=\rho(r, b)$. Let
$$
 \varepsilon=\left\{\begin{array}{ll}
2 & \text { if } r \text { and }\left\lceil\frac{r}{b}\right\rceil \text { has same parity,} \\
1 & \text { otherwise, }
\end{array} \text { and } \eta=\Big\lceil\frac{r}{b}\Big\rceil-\varepsilon\right.\text {.}
$$
Denote by $\overline{G}$ the complement graph of $G$, and $C_{n}$ the cycle on $n$ vertices, respectively. Define
$$
H_{r, \eta}= \begin{cases}\mathrm{K}_{r+1-\eta} \nabla \overline{\frac{\eta}{2} \mathrm{K}_2} & \text { if } r \text { is even, } \\ \overline{\mathrm{C}_\eta} \nabla \overline{\frac{r+2-\eta}{2} \mathrm{K}_2} & \text { if } r \text { is odd. }\end{cases}
$$
Consider the vertex partition $\{V(\overline{\mathrm{C}_\eta}), V(\overline{\frac{r+2-\eta}{2} \mathrm{K}_2})\}$ of $H_{r, \eta}$, then the quotient matrix of $A(H_{r, \eta})$ is
$$
B=\left(\begin{array}{cc}
\eta-3 & r+2-\eta \\
\eta & r-\eta
\end{array}\right).
$$
The characteristic polynomial of $B$ is
$$
p(x)=(x-\eta+3)(x-r+\eta)-(r+2-\eta) \eta.
$$
Since the vertex partition is equitable, $\rho(H_{r, \eta})=\lambda_1(B)=\frac{r-3+\sqrt{(r+3)^2-4 \eta}}{2}$. When $r$ is even, a similar result can also be obtained.

\begin{thm}[Kim, O, Park and Ree\cite{Kim-O-Park-Ree}]\label{Kim-O-Park-Ree}
Let $G$ be an $r$-regular graph of even order, where $r\geq 3$. If
 $$\lambda_{3}(G)< \rho(r,b),$$
 then $G$ has an odd $[1,b]$-factor, where
 $$
\rho(r,b)= \left\{
\begin{array}{ll}
\frac{r-2+\sqrt{(r+2)^{2}-4(\lceil\frac{r}{b}\rceil-2)}}{2}& \mbox{ if $r$ is even and $\lceil\frac{r}{b}\rceil$ is even,}\\
 \frac{r-2+\sqrt{(r+2)^{2}-4(\lceil\frac{r}{b}\rceil-1)}}{2}& \mbox{if $r$ is even and $\lceil\frac{r}{b}\rceil$ is odd,}\\
\frac{r-3+\sqrt{(r+3)^{2}-4(\lceil\frac{r}{b}\rceil-2)}}{2} & \mbox{if $r$ is odd and $\lceil\frac{r}{b}\rceil$ is odd,}\\
 \frac{r-3+\sqrt{(r+3)^{2}-4(\lceil\frac{r}{b}\rceil-1)}}{2} & \mbox{if $r$ is odd and $\lceil\frac{r}{b}\rceil$ is even}.
\end{array}
\right.
$$
\end{thm}

Motivated by the work of Amahashi\cite{A.A} and O\cite{S.O}, Li and Miao\cite{Li-Miao2022} considered the edge conditions for a connected graph to contain an odd $[1,b]$-factor.

\begin{thm}[Li and Miao\cite{Li-Miao2022}]
Let $G$ be a connected graph of even order $n$ and let $b\leq \frac{n}{2}-1$ be an odd integer.
\begin{enumerate}[(i)]
\item For $n=4$ or $n \geq 10$, if $e(G)>b+1+{n-b-1\choose 2}$, then $G$ has an odd $[1,b]$-factor.
\item For $n=6$, if $e(G)>9$, then $G$ contains a 1-factor; if $e(G)>5$, then $G$ has an odd $[1,3]$-factor.
\item For $n=8$, if $e(G)>18$, then $G$ contains a 1-factor; if $e(G)>10$, then $G$ contains an odd $[1,3]$-factor; if $e(G)>7$, then $G$ has an odd $[1,5]$-factor.
\end{enumerate}
\end{thm}

From the perspective of spectral radius condition, they also gave the following result.
\begin{thm}[Li and Miao\cite{Li-Miao2022}]\label{Li-Miao2022}
Let $G$ be a connected graph of even order $n$ and let $b\leq \frac{n}{2}-1$ be an odd integer. Assume that $\theta(n)$ is the largest root of $x^3+(b-n+3) x^2-(n-1) x-(b+1)(b-n+3)=0$.
\begin{enumerate}[(i)]
\item For $n=4$ or $n \geq 8$, if $\rho(G)>\theta(n)$, then $G$ has an odd $[1,b]$-factor.
\item For $n=6$, if $\rho(G)>\frac{1+\sqrt{33}}{2}$, then $G$ has a 1-factor; if $\rho(G)>\theta(6)$, then $G$ has an odd $[1,3]$-factor.
\end{enumerate}
\end{thm}
Observe that $K_1\nabla (K_{n-b-2}\cup (b+1)K_1)$ contains no odd $[1,b]$-factors and $K_2\nabla 4K_1$ contains no 1-factors. Moreover, $\rho(K_1\nabla (K_{n-b-2}\cup (b+1)K_1))=\theta(n)$ and $\rho(K_2\nabla 4K_1)$ $=\frac{1+\sqrt{33}}{2}$. Thus, the bounds in Theorem \ref{Li-Miao2022} are best possible.
Later, Fan, Lin and Lu\cite{Fan-Lin-Lu} generalized their result by adding the condition of minimum degree in a graph.
\begin{thm}[Fan, Lin and Lu\cite{Fan-Lin-Lu}]\label{Fan-Lin-Lu}
Suppose that $G$ is a connected graph of even order $n\geq \max\{4(b+1)\delta+4, b\delta^3+\delta\}$ with minimum degree $\delta$. If
$$\rho(G)\geq\rho(K_{\delta} \nabla (K_{n-(b+1)\delta-1}\cup (b\delta+1)K_{1})),$$
then $G$ has an odd $[1,b]$-factor unless $G\cong K_{\delta} \nabla (K_{n-(b+1)\delta-1}\cup (b\delta+1)K_{1})$.
\end{thm}

The following fundamental theorem provides a sufficient and necessary condition for the existence of a $[1,b]$-factor for $b\geq 2$. Recall that $i(H)$ is the number of isolated vertices of a graph $H$.
\begin{thm}[Berge and Las Vergnas\cite{Berge-Vergnas}]\label{Berge-Vergnas}
Let $G$ be a graph and let $b\geq 2$ be an integer. Then $G$ has a $[1,b]$-factor if and only if
$$i(G-S)\leq b|S|$$
for all subset $S\subseteq V(G)$.
\end{thm}

By Theorem \ref{Berge-Vergnas}, Fan, Lin and Lu\cite{Fan-Lin-Lu} obtained the following result.
\begin{thm}[Fan, Lin and Lu\cite{Fan-Lin-Lu}]
Suppose that $G$ is a connected graph of order $n$ with minimum degree $\delta$. If $n\geq 4(b+1)\delta+4$ and $$\rho(G)\geq\rho(K_{\delta} \nabla (K_{n-(b+1)\delta-1}\cup (b\delta+1)K_{1}))$$
 where $b\geq 2$, then $G$ has a $[1,b]$-factor unless $G\cong K_{\delta} \nabla (K_{n-(b+1)\delta-1}\cup (b\delta+1)K_{1})$.
\end{thm}

\subsection{$[a,b]$-factors}

An \textit{$[a,b]$-factor} of a graph $G$ is a spanning subgraph $H$ such that $a\leq d_{H}(v)\leq b$ for each $v\in V(G)$. A $[k,k]$-factor is called a $k$-factor. Tutte \cite{Tutte1952} obtained the well-known $k$-Factor Theorem in 1952.
\begin{thm}[Tutte\cite{Tutte1952}]
Let $k\geq 1$ be an integer and let $G$ be a graph. Then $G$ has a $k$-factor if and only if for all disjoint subsets $S, T \subseteq V(G)$,
$$
\delta_{G}(S,T)=k|S|+\sum_{x \in T}d_{G}(x)-k|T|-e_G(S, T)-q_G(S, T) \geq 0,
$$
where $q_G(S, T)$ denotes the number of components $C$ of $G-(S \cup T)$ such that $k|C|+e_G(V(C), T) \equiv 1(\bmod~2)$. Moreover, $\delta_{G}(S,T)\equiv k|V(C)|(\bmod~2)$.
\end{thm}

In 1891, Petersen\cite{Petersen1891} demonstrated that every $2r$-regular graph possesses a $2k$-factor for $0 \leq k \leq r$. One might naturally ask whether an even regular graph can have an odd factor, or an odd regular graph can have either an odd or even factor.
 Based on the Tutte's $k$-Factor Theorem, Lu\cite{Lu2012} gave an answer to this question in terms of the third largest eigenvalue.
\begin{thm}[Lu\cite{Lu2012}]
Suppose that $r$ and $k$ are two integers such that $1 \leq k<r$. Let $G$ be a connected $r$-regular graph of order $n$. Let $m$ be an integer such that $1 \leq m \leq r$ and
$$
\lambda_3(G) \leq r-\frac{m-1}{r+1}+\frac{1}{(r+1)(r+2)} .
$$
Let $m^* \in\{m, m+1\}$ such that $m^* \equiv 1(\bmod~2)$. If one of the following conditions holds, then $G$ has a $k$-factor.
\begin{enumerate}[(i)]
\item  $r$ is even, $k$ is odd, $n$ is even, and $r / m \leq k \leq r(1-1/m)$;
\item  $r$ is odd, $k$ is even, and $k \leq r(1-1 / m^*)$;
\item  both $r$ and $k$ are odd and $r / m^* \leq k$.
\end{enumerate}
\end{thm}

Using the eigenvalue interlacing Theorem, Lu\cite{Lu2010} slightly improved the above upper bounds (in terms of $r$ and $k$) on the third largest eigenvalue.
\begin{thm}[Lu\cite{Lu2010}]
Let $m$ be an integer such that $m^* \in\{m, m+1\}$ and $m^* \equiv 1(\bmod~2)$, and let $G$ be a connected $r$-regular graph. Suppose that
$$
\lambda_3(G)<\left\{
\begin{array}{ll}
\frac{r-2+\sqrt{(r+2)^{2}-4(m-1)}}{2}& \mbox{if $m$ is odd,}\\
 \frac{r-3+\sqrt{(r+3)^{2}-4(m-1)}}{2}& \mbox{if $m$ is even.}
\end{array}
\right.
$$
If one of the following conditions holds, then $G$ has a $k$-factor.
\begin{enumerate}[(i)]
\item $r$ is odd, $k$ is even and $k \leq r(1-\frac{1}{m^*})$;
\item both $r$ and $k$ are odd and $\frac{r}{m^*} \leq k$.
\end{enumerate}
\end{thm}

By considering the condition of edge-connectedness, Gallai\cite{Gallai1950} showed that an $h$-edge-connected $r$-regular graph $G$ will contain a $k$-factor depending on the parity of $r$ and $k$. In another study, Gu\cite{Gu2014} studied eigenvalue conditions of $t$-edge-connected regular graphs with $k$-factors.

\begin{thm}[Gu\cite{Gu2014}]
Suppose that $G$ is a $t$-edge-connected $r$-regular graph of order $n$. Let $k$ be an integer with $1 \leq k<d$, and let $t^{\prime} \in\{t, t+1\}$ be an even number and $t^* \in\{t, t+1\}$ be odd.

\begin{enumerate}[(i)]
\item $r$ is even, $k$ is odd and $n$ is even. Let $\hat{k}=\min \{k, r-k\}$ and $m=\lceil\frac{r}{\hat{k}}\rceil$. If $d \leq \hat{k} t^{\prime}$, or, if $d>\hat{k} t^{\prime}$ and
$$
\lambda_{\lceil\frac{2 r}{r-\hat{k} t^{\prime}}\rceil}(G)< \begin{cases}\frac{r-2+\sqrt{(r+2)^2-4(m-2)}}{2} & \text { if } $m$ \text { is even, } \\ \frac{r-2+\sqrt{(r+2)^2-4(m-1)}}{2} & \text { if } $m$ \text { is odd },\end{cases}
$$
then $G$ has a $k$-factor.
\item $r$ is odd and $k$ is even. Let $m=\left\lceil\frac{r}{r-k}\right\rceil$. If $r \leq(r-k) t^*$, or, if $r>(r-k) t^*$ and
$$
\lambda_{\lceil\frac{2 r}{r-(r-k) t^*}\rceil}(G)< \begin{cases}\frac{r-3+\sqrt{(r+3)^2-4(m-1)}}{2} & \text { if } $m$ \text { is even, } \\ \frac{r-3+\sqrt{(r+3)^2-4(m-2)}}{2} & \text { if } $m$ \text { is odd },\end{cases}
$$
then G has a $k$-factor.
\item both $r$ and $k$ are odd. Let $m=\left\lceil\frac{r}{k}\right\rceil$. If $r \leq k t^*$, or, if $r>k t^*$ and
$$
\lambda_{\lceil\frac{2 r}{r-k t^*}\rceil}(G)< \begin{cases}\frac{r-3+\sqrt{(r+3)^2-4(m-1)}}{2} & \text { if } $m$ \text { is even, } \\ \frac{r-3+\sqrt{(r+3)^2-4(m-2)}}{2} & \text { if } $m$ \text { is odd, }\end{cases}
$$
then $G$ has a $k$-factor.
\end{enumerate}
\end{thm}

Let $G$ be a graph, and $g, f$ be two integer-valued functions defined on $V(G)$. A spanning subgraph $H$ is called a \textit{$(g, f)$-factor} of $G$ if $g(v) \leq d_H(v) \leq f(v)$ for any $v \in V(G)$. In 1970, Lov\'{a}sz\cite{Lovasz1970} provided the following theorem which generalizes the criteria of other factors, such as 1-factors, $k$-factors, $[a,b]$-factors.
\begin{thm}[Lov\'{a}sz\cite{Lovasz1970}]\label{Lovasz1970}
Let $G$ be a graph and $g, f$ be integer-valued functions defined on $V(G)$ such that $g(v) \leq f(v)$ for any $v \in V(G)$. Then $G$ has a $(g, f)$-factor if and only if for all disjoint subsets $S, T \subseteq V(G)$,
$$
\sum_{s \in S} f(s)+\sum_{t \in T}\left(d_G(t)-g(t)\right)-e_G(S, T)-q_G(S, T) \geq 0,
$$
where $q_G(S, T)$ denotes the number of components $C$ of $G-(S \cup T)$ such that $g(v)=f(v)$ for any $v \in V(C)$ and $\sum_{v \in V(C)} f(v)+e_G(V(C), T) \equiv 1(\bmod~2)$.
\end{thm}

If $f(x)\equiv g(x) (\bmod~ 2)$ for all $x\in V(G)$, a $(g,f)$-factor $F$ with $d_{F}(x)\equiv f(x) (\bmod~2)$ for all $x\in V(G)$ is called a parity $(g,f)$-factor. The following criterion can be easily deduced from Lov\'{a}sz $(g,f)$-Factor Theorem.

\begin{cor}[Lov\'{a}sz $(g, f)$-parity Factor Theorem]
Let $G$ be a graph, and let $g, f$ be integer-valued functions defined on $V(G)$ such that $g(v) \leq f(v)$ for any $v \in V(G)$. Then $G$ has a $(g, f)$-factor if and only if for all disjoint subsets $S, T \subseteq V(G)$,
$$
\sum_{s \in S} f(s)+\sum_{t \in T}\left(d_G(t)-g(t)\right)-e_G(S, T)-q_G(S, T) \geq 0,
$$
where $q_{G}(S, T)$ is the number of components $C$ of $G-(S \cup T)$ such that $\sum_{v \in V(C)} f(v)+e_G(V(C), T)\equiv 1(\bmod~2)$.
\end{cor}

When $g(x)=a$ and $f(x)=b$ for all $x\in V(G)$, an even (or odd) $[a,b]$-factor of a graph $G$ is a spanning subgraph $H$ such that $d_{H}(v)$ is even (or odd) and $a\leq d_{H}(v)\leq b$ for each $v\in V(G)$. In 2022, O\cite{O2022} proved upper bounds (in terms of $a,b$ and $r$) for certain eigenvalues (in terms of $a,b,r$ and $h$) in an $h$-edge-connected $r$-regular graph $G$ to guarantee the existence of an even (or odd) $[a,b]$-factor. Denote by $r_{ab}=\min\{r-a,b\}$ and
 \begin{equation*}
 \eta=\left\{
\begin{array}{ll}
\lceil\frac{r}{r_{ab}}\rceil-1 & \mbox{if $r$ is even, $a$ and $b$ are odd, and $\lceil\frac{r}{r_{ab}}\rceil$ is odd}, \\
\lceil\frac{r}{r_{ab}}\rceil-2 & \mbox{if $r$ is even, $a$ and $b$ are odd, and $\lceil\frac{r}{r_{ab}}\rceil$ is even}, \\
\lceil\frac{r}{b}\rceil-1 & \mbox{if $r$ is odd, $a$ and $b$ are odd, and $\lceil\frac{r}{b}\rceil$ is even}, \\
\lceil\frac{r}{b}\rceil-2 & \mbox{if $r$ is odd, $a$ and $b$ are odd, and $\lceil\frac{r}{b}\rceil$ is odd}, \\
\lceil\frac{r}{r-a}\rceil-1 & \mbox{if $r$ is odd, $a$ and $b$ are odd, and $\lceil\frac{r}{r_{ab}}\rceil$ is even}, \\
\lceil\frac{r}{r-a}\rceil-2 & \mbox{if $r$ is odd, $a$ and $b$ are odd, and $\lceil\frac{r}{r_{ab}}\rceil$ is odd}.
\end{array}
\right.
\end{equation*}
Let
 \begin{equation*}
\rho(r,a,b)=\left\{
\begin{array}{ll}
\frac{r-2+\sqrt{(r+2)^2-4\eta}}{2} & \mbox{if both $r$ and $\eta$ are even}, \\
\frac{r-3+\sqrt{(r+3)^2-4\eta}}{2} & \mbox{if both $r$ and $\eta$ are odd and $\eta\geq 3$}, \\
\theta & \mbox{if $r$ is odd and $\eta=1$},
\end{array}
\right.
\end{equation*}
where $\theta$ is the largest root of $x^3-(r-2)x^2-2rx+r-1=0$.

\begin{thm}[O\cite{O2022}]\label{thm::1.4}
Let $r, a, b, h, h^{\prime}$, and $h^*$ be positive integers such that $r \geq 3, a \leq b<r$, $h \leq r$, $h^{\prime} \in\{h, h+1\}$ is an even number, and $h^* \in\{h, h+1\}$ is an odd number. Suppose that $G$ is an $h$-edge-connected $r$-regular graph of order $n$.
\begin{enumerate}[(i)]
\item  For even $r$, odd $a, b$, and even $n$, if $r \leq r_{a b} h^{\prime}$, or if $r>r_{a b} h^{\prime}$ and
$$
\lambda_{\lceil\frac{2 r}{r-r_{ab} h'}\rceil}(G)<\rho(r, a, b),
$$
then $G$ has an odd $[a, b]$-factor;
\item  For both odd $r$ and odd $a, b$, if $r \leq b h^*$, or, if $r>b h^*$ and
$$
\lambda_{\lceil\frac{2 r}{r-b h^*}\rceil}(G)<\rho(r, a, b),
$$
then $G$ has an odd $[a, b]$-factor;
\item For odd $r$ and even $a, b$, if $r \leq(r-a) h^*$, or, if $r>(r-a) h^*$ and
$$
\lambda_{\lceil\frac{2 r}{r-(r-a)h^*}\rceil}(G)<\rho(r, a, b),
$$
then $G$ has an even $[a, b]$-factor.
\end{enumerate}
\end{thm}

 Kim and O\cite{Kim-O2023} extended the result in \cite{O2022} to general factors, and investigated the sufficient conditions for a graph to have a $(g, f)$-parity factor in terms of the minimum degree, edge-connectivity and eigenvalues.
For positive integers $r$ and $h$ with $r>h$, let
$$
\rho(r, h)= \begin{cases}\frac{r-2+\sqrt{(r+2)^2-8\lfloor h / 2\rfloor}}{2} & \text { if } r \text { or } h \text { is even, } \\ \frac{r-3+\sqrt{(r+3)^2-4 h}}{2} & \text { if both } r, h \text { are odd and } h \geq 3, \\ \mu & \text { if } r \text { is odd and } h=1,\end{cases}
$$
where $\mu$ is the largest real root of $x^3-(r-2) x^2-2 r x+r-1=0$.

\begin{thm}[Kim and O\cite{Kim-O2023}]\label{Kim-O2023}
Suppose that $G$ is an $h$-edge-connected graph and $\theta$ is a real number with $0<\theta<1$. Let $g$ and $f$ be integer-valued functions on $V(G)$ such that $g(v) \leq \theta d_G(v) \leq f(v)$ and $g(v) \equiv f(v)~(\bmod~2)$ for all $v \in V(G)$ and $\sum_{v \in V(G)} f(v)$ is even. Let $\theta^*=\min \{\theta, 1-\theta\}$ and let $h_e$ and $h_o$ be even and odd integers in $\{h, h+1\}$, respectively. If one of (a)-(e) holds, then $G$ has a $(g, f)$- parity factor.
\begin{enumerate}[(a)]
\item i. $h \geq 1 / \theta^*$, or\\
ii. $h<1 / \theta^* \leq \delta(G)$ and $\lambda_{\lceil\frac{2}{1-\theta^* h}\rceil}<\rho(\delta(G),\lceil 1 / \theta^*\rceil-1)$.
\item $d_G(v)$ and $f(v)$ are even for all $v \in V(G)$.
\item $d_G(v)$ is even for all $v \in V(G)$, and one of the following holds.\\
i. $h_e \geq 1 / \theta^*$.\\
ii. $h_e<1 / \theta^* \leq \delta(G)$, and $\lambda_{\lceil\frac{2}{1-\theta^* h_e}\rceil}<\rho(\delta(G),\lceil 1 / \theta^*\rceil-1)$.
\item $f(v)$ is even for all $v \in V(G)$, and one of the following holds.\\
i. $h_o \geq 1 /(1-\theta)$.\\
ii. $h_o<1 /(1-\theta) \leq \delta(G)$ and $\lambda_{\lceil\frac{2}{1-(1-\theta) h_o}\rceil}<\rho(\delta(G),\lceil 1 /(1-\theta)\rceil-1)$.
\item $d_G(v) \equiv f(v)(\bmod~2)$ for all $v \in V(G)$, and one of the following holds.\\
i. $h_o \geq 1 / \theta$.\\
ii. $h_o<1 / \theta \leq \delta(G)$ and $\lambda_{\lceil\frac{2}{1-\theta h_o}\rceil}<\rho(\delta(G),\lceil 1 / \theta\rceil-1)$.
\end{enumerate}
\end{thm}

For $n \geq 2 x>0$, Cho, Hyun, O and Park\cite{Cho-Hyun-O-Park} proved that a complete bipartite graph $K_{x, n-x}$ has an $[a, b]$-factor if and only if
$$
\rho(K_{x, n-x}) \geq\left\{\begin{array}{cl}
\sqrt{a(n-a)} & \text { if } n<a+b, \\
\sqrt{a b}\cdot \frac{n}{a+b} & \text { if } n \geq a+b .
\end{array}\right.
$$
Among graphs of order $n$ without an $[a, b]$-factor, they guessed that the graph $K_{a-1}\nabla(K_{1}\cup K_{n-a})$ attains the maximum spectral radius. Note that there are $n-a$ vertices with degree $n-2$, $a-1$ vertices with degree $n-1$, and one vertex with degree $a-1$ in $K_{a-1}\nabla(K_{1}\cup K_{n-a})$. Thus, $K_{a-1}\nabla(K_{1}\cup K_{n-a})$ contains no $[a, b]$-factors. From this, they posed a conjecture regarding the spectral radius condition for the existence of an $[a,b]$-factor in graphs as follows.
\begin{conj}[Cho, Hyun, O and Park \cite{Cho-Hyun-O-Park}]\label{Cho-Hyun-O-Park}
Let $a\cdot n$ be an even integer at least 2 where $n\geq a+1$, and let $G$ be a graph of order $n$. If
$$\rho(G)>\rho(K_{a-1}\nabla(K_{1}\cup K_{n-a})),$$
then $G$ has an $[a,b]$-factor.
\end{conj}

Li and Cai \cite{Li-Cai} showed that a graph $G$ of order $n \geq 2 a+b+\frac{a^2-a}{b}$ contains an $[a, b]$-factor if the maximum degree of any two non-adjacent vertices of $G$ is greater than $\frac{a n}{a+b}$. In light of this condition, Fan, Lin and Lu\cite{Fan-Lin-Lu} confirmed Conjecture \ref{Cho-Hyun-O-Park} for $n\geq 3a+b-1$.

\begin{thm}[Fan, Lin and Lu\cite{Fan-Lin-Lu}]
 Let $a\cdot n$ be an even integer, where $n\geq 3a+b-1$ and $b\geq a\geq 1$, and let $G$ be a graph of order $n$. If $$\rho(G)>\rho(K_{a-1}\nabla(K_{1}\cup K_{n-a})),$$
then $G$ has an $[a,b]$-factor.
\end{thm}

For nonnegative integers $a$ and $b$, an $[a,b]$-factor is a $(g,f)$-factor with $g\equiv a$ and $f\equiv b$. By using Theorem \ref{Lovasz1970} as a tool, Wei and Zhang\cite{Wei-Zhang} confirmed the Conjecture \ref{Cho-Hyun-O-Park} completely.


\section{Eigenvalues and component factors}
In this section, we study  component factors. 
For a family $\mathcal{S}$ of graphs, a $\mathcal{S}$-factor of a graph $G$ is a spanning subgraph of $G$ such that each of its components is isomorphic to an element of $\mathcal{S}$.
\subsection{Path factors}
A path factor of a graph $G$ is a spanning subgraph of $G$ in which each component
is a path of order at least $2$. Since every path of order at least four has a $\{P_2,P_3\}$-factor, a graph $G$ has a path factor if and only if $G$ has a $\{P_2,P_3\}$-factor. Akiyama, Jin and Era\cite{Akiyama-Jin-Era} presented the characterization of a path factor.
\begin{thm}[Akiyama, Jin and Era\cite{Akiyama-Jin-Era}]\label{Akiyama-Jin-Era}
A graph $G$ has a $\{P_2,P_3\}$-factor if and only if $$i(G-S) \leq 2|S|$$ for all $S \subseteq V(G)$.
\end{thm}
Li and Miao\cite{Li-Miao} used Theorem \ref{Akiyama-Jin-Era} to provide a spectral radius condition that ensures a graph $G$ contains a $\{P_2,P_3\}$-factor.
\begin{thm}[Li and Miao \cite{Li-Miao}]\label{Li-Miao}
Let $G$ be a connected graph of order $n$.
\begin{enumerate}[(i)]
\item For $n \geq 4$ and $n \neq 7$, if $\rho(G)>\theta(n)$, then $G$ has a $\{P_2,P_3\}$-factor, where $\theta(n)$ is the largest root of $x^3-(n-5)x^2+(1-n) x+3(n-5)=0$.
\item For $n=7$, if $\rho(G)>\frac{1+\sqrt{41}}{2}$, then $G$ has a $\{P_2,P_3\}$-factor.
\end{enumerate}
\end{thm}
Note that $\rho({K}_1\nabla(K_{n-4}\cup 3K_1))=\theta(n)$ and $\rho({K}_2 \nabla 5K_1)=\frac{1+\sqrt{41}}{2}$. The conditions in Theorem \ref{Li-Miao} are sharp since none of ${K}_1\nabla(K_{n-4}\cup 3K_1)$ and $K_2\nabla 5K_1$ contains no $\{P_2,P_3\}$-factors. Subsequently, Zhang\cite{W.Zhang} improved the result presented in \cite{Li-Miao} by incorporating the minimum degree condition.

\begin{thm}[Zhang \cite{W.Zhang}]
Suppose that $n \geq 2$ and $t \geq 1$ are two integers. Let $G$ be a connected graph of order $n$ with minimum degree $\delta(G) \geq t$ and without a $\{P_2,P_3\}$-factor. Then $n \geq 3 t+1$ and
$$
\rho(G) \leq \max\Big\{\rho(K_t \nabla(K_{n-1-3 t} \cup \overline{K_{2t+1}})), \rho\Big(K_{\lfloor\frac{n-1}{3}\rfloor} \nabla\Big(K_{n-1-3\lfloor\frac{n-1}{3}\rfloor} \cup \overline{K_{2\lfloor\frac{n-1}{3}\rfloor+1}}\Big)\Big)\Big\},
$$
with equality if and only if
$G\cong K_t \nabla(K_{n-1-3 t}\cup\overline{K_{2t+1}})$ or $G\cong K_{\lfloor\frac{n-1}{3}\rfloor} \nabla\Big(K_{n-1-3\lfloor\frac{n-1}{3}\rfloor} \cup $ $\overline{K_{2\lfloor\frac{n-1}{3}\rfloor+1}}\Big)\Big.$.
\end{thm}

A graph $G$ is said to be \textit{factor-critical} if $G-v$ has a perfect matching for every vertex $v \in V(G)$. Let $G$ be a factor-critical graph of order $n\geq 3$ and let $V(G)=\left\{v_1, v_2, \ldots, v_n\right\}$. Add $n$ new vertices $\left\{w_1, w_2, \ldots, w_n\right\}$ to $G$ together with edges $v_i w_i$ for $1 \leq i \leq n$. The resulting graph on $2n$ vertices is called a sun. Let $sun(G)$ denote the number of components of $G$ which are sun graphs.
\begin{thm}[Kaneko\cite{Kaneko2003}, Kano, Katona and Kir\'{a}ly\cite{Kano-Katona-Kiraly}]\label{Kano-Katona-Kiraly-Kaneko}
A graph $G$ has a $\{P_3,P_4,P_5\}$-factor if and only if
$$sun(G-S) \leq 2|S|$$
 for any $S \subseteq V(G)$.
\end{thm}
By employing Theorem \ref{Kano-Katona-Kiraly-Kaneko} as a valuable tool, Zhang\cite{W.Zhang} gave the maximum spectral radius among all connected graph without a $\{P_3,P_4,P_5\}$-factor.
\begin{thm}[Zhang \cite{W.Zhang}]
Suppose that $n \geq 3$ and $t \geq 1$ are two integers. Let $G$ be a connected graph of order $n$ with minimum degree $\delta(G) \geq t$ and without a $\{P_3,P_4,P_5\}$-factor. Then $n\geq 3t+1$ and the following conclusions hold.
\begin{enumerate}[(i)]
\item If $n \equiv 0(\bmod~3)$, then
$$\rho(G) \leq \max\Big\{\rho(K_t \nabla (K_{n-1-3 t} \cup \overline{K_{2t+1}})), \rho\Big(K_{\frac{n-3}{3}} \nabla\Big(K_2 \cup K_2\cup \overline{K_{\frac{2n-9}{3}}}\Big)\Big)\Big\},$$
with equality if and only if $G\cong K_t \nabla(K_{n-1-3 t} \cup \overline{K_{2t+1}})$ or $G\cong K_{\frac{n-3}{3}} \nabla(K_2 \cup K_2 \cup \overline{K_{\frac{2n-9}{3}}}).$
\item If $n \equiv 2(\bmod~3)$, then
$$
\rho(G) \leq \max\Big\{\rho(K_t \nabla(K_{n-1-3 t} \cup  \overline{K_{2t+1}})), \rho\Big(K_{\frac{n-2}{3}}\nabla\Big(K_2 \cup \overline{K_{\frac{2(n-2)}{3}}}\Big)\Big)\Big\},
$$
with equality if and only if $G\cong K_t \nabla(K_{n-1-3 t} \cup  \overline{K_{2t+1}})$ or $G\cong K_{\frac{n-2}{3}}\nabla(K_2 \cup \overline{K_{\frac{2(n-2)}{3}}})$.
\item If $n \equiv 1(\bmod~3)$, then
$$
\rho(G) \leq \max\Big\{\rho(K_t\nabla(K_{n-1-3t} \cup \overline{K_{2t+1}})), \rho\Big(K_{\frac{n-1}{3}}\nabla \overline{K_{\frac{2n+1}{3}}}\Big)\Big\},
$$
with equality if and only if $G\cong K_t\nabla(K_{n-1-3t} \cup \overline{K_{2t+1}})$ or $G\cong K_{\frac{n-1}{3}}\nabla\overline{K_{\frac{2n+1}{3}}}$.
\end{enumerate}
\end{thm}

\subsection{Star factors}
Let $G$ be a graph and let $k \geq 2$ be an integer. A \textit{ $\{K_{1, j}: 1 \leq j \leq k\}$-factor} of a graph $G$ is a spanning subgraph of $G$ such that each component is isomorphic to a member in $\{K_{1, j}: 1 \leq j \leq k\}$. Miao and Li\cite{Miao-Li} established a lower bound on the spectral radius conditions to ensure the existence of a $\{K_{1, j}: 1 \leq j \leq k\}$-factor in a graph. Furthermore, they also constructed certain extremal graphs to demonstrate that all bounds presented below are indeed best possible.

\begin{thm}[Miao and Li \cite{Miao-Li}]
Let $k \geq 2$ and $n \geq k+2$ be two integers, and let $G$ be a connected graph of order $n$. Assume that $\theta(n, k)$ is the largest root of $x^3-(n-k-3) x^2-(n-1) x-(k+3-n)(k+1)=0$.
\begin{enumerate}[(i)]
\item  For $k=2$, if $\rho(G)>\theta(n, 2)$ with $n \geq 4$ and $n \neq 7$, then $G$ has a $\left\{K_{1,1}, K_{1,2}\right\}$-factor; if $n=7$ and $\rho(G)>\frac{1+\sqrt{41}}{2}$, then $G$ has $a\left\{K_{1,1}, K_{1,2}\right\}$-factor.
\item  For $k=3$, if $\rho(G)>\theta(n, 3)$ with $n \geq 5$ and $n \neq 9$, then $G$ has a $\left\{K_{1,1}, K_{1,2}, K_{1,3}\right\}$-factor; if $n=9$ and $\rho(G)>\frac{1+\sqrt{57}}{2}$, then $G$ has $a\left\{K_{1,1}, K_{1,2}, K_{1,3}\right\}$-factor.
\item  For $k \geq 4$, if $\rho(G)>\theta(n, k)$, then $G$ has a $\left\{K_{1, j}: 1 \leq j \leq k\right\}$-factor.
\end{enumerate}
\end{thm}

\subsection{Cycle factors}
For $k\geq 3$, Luo, Liu and Ao\cite{Luo-Liu-Ao} proved a tight spectral condition to guarantee the existence of a $\left\{K_2,\left\{C_k\right\}\right\}$-factor in a graph with minimum degree.

\begin{thm}[Luo, Liu and Ao\cite{Luo-Liu-Ao}]
Let $k \geq 3$ be an integer, and let $G$ be a graph of order $n \geq 5 \delta+6$ with minimum degree $\delta \geq 1$. If
$$
\rho(G) \geq \rho\left(K_\delta \nabla\left(K_{n-2 \delta-1}\cup(\delta+1) K_1\right)\right),
$$
then $G$ has a $\left\{K_2,\left\{C_k\right\}\right\}$-factor unless $G \cong K_\delta \nabla\left(K_{n-2 \delta-1}\cup(\delta+1) K_1\right)$.
\end{thm}

\subsection{Connected factors}

 A graph $G$ has a $(g,f)$-factor,  but it is not necessarily the case that $G$ has a connected $(g,f)$-factor. The problem of determining whether a graph $G$ has a connected $(g,f)$-factor is NP-complete\cite{Garey-Johnson}. As is well known, it remains so in the case of $g\equiv f\equiv 2$. A cycle $C$ of a graph $G$ is a \textit{Hamiltonian cycle} if it passes once and only once through every vertex of $G$; likewise, a path is a \textit{Hamiltonian path} if it passes every vertex exactly once.  A graph is called a Hamiltonian graph if it contains a Hamiltonian cycle. Obviously, A Hamiltonian cycle is a connected 2-factor and a Hamiltonian path is a connected $[1,2]$-factor, respectively. Hence, the factor theory of graphs can be viewed as an extension of the Hamiltonian cycle problem.

The following well-known theorem was given by Ore, which is an edge condition for the existence of a Hamiltonian cycle in a graph.
\begin{thm}[Ore \cite{Ore}]\label{Ore}
Let $G$ be a graph of order $n$. If
$$e(G)>{n-1\choose2}+1.$$
then $G$ has a Hamiltonian cycle.
\end{thm}

In 2010, Fiedler and Nikiforov\cite{Fiedler-Nikiforov} gave spectral radius conditions for the existence of a Hamiltonian cycle or Hamiltonian path in a graph, respectively.
\begin{thm}[Fiedler and Nikiforov \cite{Fiedler-Nikiforov}]
Let $G$ be a graph of order $n$. If
$$\rho(G)\geq n-2,$$
 then $G$ has a Hamiltonian path unless $G\cong K_{n-1}\cup K_1$. If strict inequality holds, then $G$ has a Hamiltonian cycle unless $G\cong K_{n-1}+e$, where $K_{n-1}+e$ denotes the complete graph on $n-1$ vertices with a pendent edge.
\end{thm}

Note that $\delta(G)\geq 2$ is a trivial necessary condition for a graph to be Hamiltonian.
\begin{thm}[Ning and Ge \cite{Ning-Ge}]\label{Ning-Ge}
Let $G$ be a connected graph of order $n\geq 14$ with minimum degree $\delta(G)\geq 2$. If
\begin{equation*}
\begin{aligned}
\rho(G)\geq \rho(K_2\nabla(K_{n-4}\cup 2K_1)),
\end{aligned}
\end{equation*}
then $G$ has a Hamiltonian cycle unless $G\cong K_2\nabla(K_{n-4}\cup 2K_1)$.
\end{thm}
For $n\leq 13$, $K_2\nabla(K_{n-4}\cup 2K_1)$ is not necessarily the extremal graph that corresponds to the lower bound. For example, when $n=7$, $\rho(K_3 \nabla 4 K_1)>\rho(K_2\nabla(K_3\cup 2K_1)$ and $K_3 \nabla 4 K_1$ contains no Hamiltonian cycles, and when $n=9$, $\rho(K_4 \nabla 5 K_1)>\rho(K_2\nabla(K_5\cup 2K_1))$ and $K_4 \nabla 5K_1$ contains no Hamiltonian cycles. Moreover, Ning and Ge\cite{Ning-Ge} proposed the conjecture that the condition $n\geq 14$ in Theorem \ref{Ning-Ge} could be improved to $n\geq 10$, which was later verified by Chen, Hou and Qian\cite{Chen-Hou-Qian}.
In the same paper, they also considered a spectral radius condition for the existence of a Hamiltonian path in graphs.
\begin{thm}[Ning and Ge \cite{Ning-Ge}]
Let $G$ be a graph of order $n \geq 4$ with minimum degree $\delta(G) \geq 1$. If
$$\rho(G)>n-3,$$
then $G$ has a Hamiltonian path unless $G\in \{K_1\nabla(K_{n-3}\cup 2 K_1), K_2\nabla 4 K_1, K_1\nabla$ $ (K_{1,3}\cup K_1)\}$.
\end{thm}

Soon after, Benediktovich\cite{Benediktovich} improved the results in \cite{Fiedler-Nikiforov,Ning-Ge}.
\begin{thm}[Benediktovich \cite{Benediktovich}]\label{Benediktovich}
Let $G$ be a graph of order $n\geq 9$ with minimum degree $\delta(G) \geq 2$. If
$$\rho(G)>n-3,$$
then $G$ has a Hamiltonian cycle unless $G\in \{K_4\nabla 5K_1, K_3\nabla(K_{1,4}\cup K_1), K_1\nabla (K_{n-3}\cup K_2), K_2\nabla(K_{n-4}\cup 2K_1)\}$.
\end{thm}

The celebrated Dirac theorem\cite{Dirac1952} states that every graph of order $n\geq 3$ with minimum degree at least $n/2$ has a Hamiltonian cycle. By taking the minimum degree as a parameter, Erd\H{o}s\cite{Erdos} gave a sufficient condition for a graph to contain a Hamiltonian cycle which generalized Ore's theorem\cite{Ore}.
\begin{thm}[Erd\H{o}s \cite{Erdos}]\label{Erdos}
Let $G$ be a graph of order $n\geq 14$ with minimum degree $\delta(G)\geq k$ where $1\leq k\leq (n-1)/2$. If
$$e(G)>\max\Big\{{n-k\choose2}+k^2, {\lceil\frac{n+1}{2}\rceil\choose2}+\Big\lfloor\frac{n-1}{2}\Big\rfloor^2\Big\},$$
then $G$ has a Hamiltonian cycle.
\end{thm}

 Li and Ning\cite{Li-Ning} presented the spectral analogues of Erd\H{o}s theorem.
\begin{thm}[Li and Ning \cite{Li-Ning}]\label{Li-Ning-1}
Let $k$ be an integer, and let $G$ be a graph of order $n$ with minimum degree $\delta(G) \geq k$.
\begin{enumerate}[(i)]
\item If $k \geq 1$, $n\geq\max\{6k+5,(k^2+6k+4)/2\}$ and $$\rho(G)\geq \rho(K_{k}\nabla(K_{n-2k}\cup kK_1)),$$
then $G$ has a Hamiltonian cycle unless $G\cong K_{k}\nabla(K_{n-2k}\cup kK_1)$.
\item If  $k \geq 0$, $n\geq\max\{6k+10,(k^2+7k+8)/2\}$ and $$\rho(G)\geq \rho(K_{k}\nabla(K_{n-2k-1}\cup (k+1)K_1)),$$
then $G$ has a Hamiltonian path unless $G\cong K_{k}\nabla(K_{n-2k-1}\cup (k+1)K_1)$.
\end{enumerate}
\end{thm}

Nikiforov\cite{Nikiforov2016} proposed the following theorem, which generalizes the Theorems \ref{Ning-Ge}
and \ref{Benediktovich}, and strengthen the Theorem \ref{Li-Ning-1} for $n$ sufficiently large.
\begin{thm}[Nikiforov \cite{Nikiforov2016}]
Let $k\geq 1$ be an integer, and let $G$ be a connected graph of order $n$ with minimum degree $\delta(G)\geq k$.
\begin{enumerate}[(i)]
\item
If $n\geq k^3+k+4$ and
$$\rho(G)\geq n-k-1,$$
then $G$ has a Hamiltonian cycle unless $G\cong K_1\nabla(K_{n-k-1}\cup K_k)$ or $G\cong K_k\nabla(K_{n-2k}\cup kK_1)$.
\item If $n \geq k^3+k^2+2 k+5$ and
$$\rho(G)\geq n-k-2,$$
then $G$ has a Hamiltonian path unless $G\cong K_{k}\nabla(K_{n-2k-1}\cup (k+1)K_1)$ or $G\cong K_{n-k-1}\cup K_{k+1}$.
\end{enumerate}
\end{thm}


A graph $G$ is \textit{$t$-tough} if $|S|\geq tc(G-S)$ for every subset $S\subseteq V(G)$ with $c(G-S)>1$, where $c(G)$ is the number of components of a graph $G$. Being 1-tough is an obvious necessary condition for a graph to be Hamiltonian\cite{Bondy-Murty}. Conversely, Chv\'{a}tal\cite{V.C} proposed the following conjecture.
\begin{conj}[Chv\'{a}tal \cite{V.C}]
There exists a finite constant $t_{0}$ such that every $t_0$-tough graph is Hamiltonian.
\end{conj}
Around this conjecture, Enomoto, Jackson, Katerinis and Saito \cite{H.E-1} obtained
that every 2-tough graph contains a 2-factor and there exist $(2-\varepsilon)$-tough graphs without a 2-factor, where $\varepsilon>0$ is an arbitrary small number, and hence without a Hamiltonian cycle. They conjectured that the value of $t_{0}$ might be 2. In 2000, this result was disproved by Bauer, Broersma and Veldman \cite{D.B-3}, and they observed that if such a $t_0$ exists, then it must be at least $\frac{9}{4}$. In general, the conjecture still remains open. By incorporating toughness and spectral conditions, Fan, Lin and Lu\cite{Fan-Lin-Lu2023} considered Chv\'{a}tal's conjecture from another perspective and determined the spectral condition to guarantee the existence of a Hamiltonian cycle among $1$-tough graphs. Let $K_{n-4}^{+3}$ be the graph obtained from $3K_{1}\cup K_{n-4}$ by adding three independent edges between $3K_{1}$ and $K_{n-4}$, and let $M_{n}=K_{1}\nabla K_{n-4}^{+3}$.

\begin{thm}[Fan, Lin and Lu\cite{Fan-Lin-Lu2023}]
Suppose that $G$ is a connected 1-tough graph of order $n\geq 18$ with minimum degree $\delta(G)\geq 2$. If $\rho(G)\geq\rho(M_{n})$, then $G$ has a Hamiltonian cycle, unless $G\cong M_{n}$.
\end{thm}

Moon and Moser\cite{Moon-Moser} pointed out that $\delta(G)>\frac{n}{2}$ is a necessary condition for a balanced bipartite on $2n$ vertices to be Hamiltonian. Similar to the work of Erd\H{o}s in Theorem \ref{Erdos}, they considered the corresponding result for the balanced bipartite graphs with minimum degree at most $\frac{n}{2}$.

\begin{thm}[Moon and Moser\cite{Moon-Moser}] \label{Moon-Moser}
Let $G$ be a balanced bipartite graph of order $2n$. If the minimum degree $\delta(G) \geq k$ for  $1 \leq k \leq \frac{n}{2}$ and
$$
e(G)>n(n-k)+k^2,
$$
then $G$ has a Hamiltonian cycle.
\end{thm}

Let $B_{n, k}$ be the bipartite graph obtained from the complete bipartite graph $K_{n, n}$ by deleting all edges in its one subgraph $K_{k, n-k}$. One can verify that the condition in Theorem \ref{Moon-Moser} is best possible since $e\left(B_{n, k}\right)=n(n-k)+k^2$ and $B_{n, k}$ contains no Hamiltonian cycles.

Liu, Shiu and Xue\cite{Liu-Shiu-Xue} gave a sufficient condition on the spectral radius for a bipartite graph to be Hamiltonian. Let $K_{n,n-2}+4e$ be a bipartite graph obtained from $K_{n, n-2}$ by adding two vertices which are adjacent to two common vertices with degree $n-2$ in $K_{n, n-2}$.

\begin{thm}[Liu, Shiu and Xue\cite{Liu-Shiu-Xue}]
Let $G$ be a balanced bipartite graph of order $2n \geq 8$ with minimum degree $\delta(G) \geq 2$. If
$$
\rho(G) \geq \sqrt{n(n-2)+4},
$$
then $G$ has a Hamiltonian cycle unless $G\cong K_{n, n-2}+4 e$.
\end{thm}

Li and Ning\cite{Li-Ning} focused on a spectral analogue about the result in Theorem \ref{Moon-Moser}.

\begin{thm}[Li and Ning\cite{Li-Ning}]\label{Li-Ning}
Suppose that $k \geq 1$ and $n \geq(k+1)^2$. Let $G$ be a balanced bipartite graph of order $2n$ with minimum degree $\delta(G) \geq k$. If
$$
\rho(G) \geq \rho(B_{n, k}),
$$
then $G$ has a Hamiltonian cycle unless $G\cong B_{n, k}$.
\end{thm}

Since $K_{n, n-k}$ is a proper subgraph of $B_{n, k}$,  $\rho(B_{n, k})>\rho(K_{n, n-k})=\sqrt{n(n-k)}$. By using this result, Ge and Ning\cite{Ge-Ning2020} obtained an improved boundary on Theorem \ref{Li-Ning}.

\begin{thm}[Ge and Ning\cite{Ge-Ning2020}]\label{Ge-Ning2020}
 Suppose that $k \geq 1$ and $n \geq k^3+2 k+4$. Let $G$ be a balanced bipartite graph of order $2n$ with minimum degree $\delta(G) \geq k$. If
$$
\rho(G) \geq \sqrt{n(n-k)},
$$
then $G$ has a Hamiltonian cycle unless $G\cong B_{n, k}$.
\end{thm}

Notably, Jiang, Yu and Fang\cite{Jiang-Yu-Fang} independently derived the result stated in Theorem \ref{Ge-Ning2020} to be valid under a relaxed condition of $n \geq \max \left\{\frac{k^3}{2}\!+\!k\!+\!3,(k+1)^2\right\}$.

A basic result in graph theory asserts that every connected graph has a spanning tree. On the other hand, a spanning tree of a graph on $n$ vertices is a connected $[1,n-1]$-factor. Many researchers investigated the existence of spanning trees under given conditions. A \textit{$k$-tree} is a spanning tree with every vertex of degree at most $k$, where $k\geq 2$ is an integer. It is not difficult to see that a $k$-tree is a connected $[1,k]$-factor. In 1989, Win \cite {S.Win} proved the following result, which gives a sufficient condition for the existence of a $k$-tree in a connected graph, see \cite{M.Ellingham} for a short proof. Let $c(G)$ be the number of components of a graph $G$.
\begin{thm}[Win \cite{S.Win}]\label{S.Win}
Let $G$ be a connected graph. If $k\geq 2$ and
$$c(G-S)\leq(k-2)|S|+2$$
for $S\subseteq V(G)$, then $G$ has a $k$-tree.
\end{thm}
Since a $2$-tree is just a Hamiltonian path, the results in\cite{Fiedler-Nikiforov,Ning-Ge,Li-Ning} provide the spectral conditions when $k=2$. For $k \geq 3$, Wong\cite{Wong2013} employed the result of Win and established the following theorem regarding the existence of $k$-trees in regular graphs.
\begin{thm}[Wong \cite{Wong2013}]\label{Wong2013}
Let $k \geq 3$ and let $G$ be a connected $r$-regular graph. If
$$\lambda_4(G) <r-\frac{r}{(k-2)(r+1)},$$
then $G$ has a $k$-tree.
\end{thm}

Cioab\u{a} and Gu\cite{Cioaba-Gu} generalized Theorem \ref{Wong2013} by taking the edge connectivity into account.
\begin{thm}[Cioab\u{a} and Gu\cite{Cioaba-Gu}]
Let $k \geq 3$ and let $G$ be a connected $r$-regular graph with edge connectivity $\kappa^{\prime}$. Let $l=r-(k-2) \kappa^{\prime}$. Each of the following statements holds.
\begin{enumerate}[(i)]
\item If $l \leq 0$, then $G$ has a $k$-tree.
\item If $l>0$ and $\lambda_{\lceil\frac{3 r}{l}\rceil}<r-\frac{r}{(k-2)(r+1)}$, then $G$ has a $k$-tree.
\end{enumerate}
\end{thm}
 Fan, Goryainov, Huang and Lin\cite{D.F} gave a spectral radius condition for the existence of a $k$-tree in a connected graph of order $n$.
\begin{thm}[Fan, Goryainov, Huang and Lin \cite{D.F}]
Let $k\geq 3$, and let  $G$ be a connected graph of order $n\geq 2k+16$. If
$$\rho(G)\geq \rho(K_{1} \nabla (K_{n-k-1}\cup k K_{1})),$$
then $G$ has a $k$-tree unless $G\cong K_{1} \nabla (K_{n-k-1}\cup kK_{1})$.
\end{thm}

For any integer $k\geq 2$, a \textit{spanning $k$-ended-tree} of a connected graph $G$ is a spanning tree with at most $k$ leaves. For an integer $l \geq 0$, the $l$-closure of a graph $G$ is the graph obtained from $G$ by successively joining pairs of nonadjacent vertices whose degree sum is at least $l$ until no such pair exists. Let $C_l(G)$ denote the $l$-closure of $G$. The following closure theorem regarding the existence of a spanning $k$-ended-tree was established by Broersma and Tuinstra.

\begin{thm}[Broersma and Tuinstra\cite{Broersma-Tuinstra}]
 Let $G$ be a connected graph of order $n$, and let $k$ be an integer with $2 \leq k \leq n-1$. Then $G$ has a spanning $k$-ended-tree if and only if the $(n-1)$-closure $C_{n-1}(G)$ of $G$ has a spanning $k$-ended-tree.
\end{thm}
Ao, Liu and Yuan\cite{Ao-Liu-Yuan} gave the following result.
\begin{thm}[Ao, Liu and Yuan\cite{Ao-Liu-Yuan}]
 Let $G$ be a connected graph of order $n$ and $k \geq 2$ be an integer. If $n \geq \max \left\{6 k+5, k^2+\frac{3}{2} k+2\right\}$ and
 $$\rho(G) \geq \rho\left(K_1 \nabla(K_{n-k-1}\cup k K_1)\right),$$
 then $G$ has a spanning $k$-ended-tree unless $G \cong K_1 \nabla(K_{n-k-1}\cup k K_1)$.
\end{thm}

Let $T$ be a spanning tree of a connected graph $G$. The leaf degree of a vertex $v \in V(T)$ is the number of leaves adjacent to $v$ in $T$.  Additionally, the leaf degree of $T$ refers to the maximum leaf degree among all the vertices of $T$. Kaneko\cite{Kaneko} gave a characterization of trees with a fixed leaf degree using the number of isolated vertices.
\begin{thm}[Kaneko\cite{Kaneko}]
Let $k \geq 1$ be an integer, and let $G$ be a connected graph. Then $G$ has a spanning tree with leaf degree at most $k$ if and only if
$$i(G-S)<(k+1)|S|$$
for any nonempty subset $S \subseteq V(G)$.
\end{thm}

Motivated by Kaneko's theorem, Ao, Liu and Yuan\cite{Ao-Liu-Yuan} presented a tight spectral condition for the existence of a spanning tree with leaf degree at most $k$ in a connected graph.

\begin{thm}[Ao, Liu and Yuan\cite{Ao-Liu-Yuan}]
Let $k \geq 1$ be an integer, and let $G$ be a connected graph of order $n \geq 2 k+12$. If
$$\rho(G) \geq \rho(K_1 \nabla(K_{n-k-2}\cup (k+1) K_1)),$$
then $G$ has a spanning tree with leaf degree at most $k$ unless $G \cong K_1 \nabla(K_{n-k-2}\cup (k+1) K_1)$.
\end{thm}

\section{Eigenvalues and factor packing problem}

As an extension of the factor existence problem, many researchers have explored the maximum number of edge-disjoint 1-factor\cite{Faudree-Gould}, Hamiltonian cycle\cite{Li2000} and spanning tree\cite{Zhou-Bu-Lai} in a graph by using graph parameters. However, there are relatively few results that investigate this problem from the perspective of eigenvalues.

\subsection{Perfect matchings}

One problem concerned perfect matching which has attracted considerable interest is that of determining the structure of graphs with a unique perfect matching\cite{X.W,D.B,S.P}. In 1985, Godsil\cite{Godsil1985} showed that the path attains the minimum smallest positive eigenvalue among all trees with a unique perfect matching, and posted an open problem to characterize all bipartite graphs with a unique perfect matching whose adjacency matrices have inverses diagonally similar to non-negative matrices. This problem was settled by Yang and Ye\cite{Yang-Ye}.
 For a unicyclic graph $G$, if $G$ contains a unique perfect matching, then $A(G)^{-1}$ is an integer matrix\cite{Akbari-Kirkland}. Let $A(G)^{-1}=\left[b_{i j}\right]$, and let $G_{+}$ be the graph associated to the matrix $A(G)^{-1}$ such that $V(G_{+})=V(G)$ and the vertices $i$ and $j$ are adjacent in $G_{+}$ if and only if $b_{i j} \neq 0$. This implies that $\tau(G)=\frac{1}{\rho\left(G_{+}\right)}$, where $\tau(G)$ denotes the smallest positive eigenvalue of $A(G)$. By using the concept of inverse graphs as a tool, Barik and Behera\cite{Barik-Behera} characterized the unique extremal graph attains the minimum $\tau$ value among all connected bipartite unicyclic graphs on $n=2m$ vertices.

 Let $\mathcal{U}_{2 m}$ be the class of connected bipartite unicyclic
graphs on $n=2m$ vertices. Consider a path $P_{2 m}=v_1v_2\cdots v_{m}v_{m+1}\cdots v_{2m}$. If $m$ is even, let $U_e$ be the graph obtained from $P_{2 m}$ by adding an edge between the vertices $v_{m-2}$ and $v_{m+3}$. If $m$ is odd, let $U_o$ be the graph obtained from the path $P_{2 m}$ by adding an edge between the vertices $v_{m-3}$ and $v_{m+2}$. It is worth mentioning that both $U_e$ and $U_o$ are very close to $P_{2 m}$.
\begin{thm}[Barik and Behera \cite{Barik-Behera}]
Let $m \geq 4$ and $U \in \mathcal{U}_{2 m}$. Then the following results hold.
\begin{enumerate}[(i)]
\item If $m$ is even, then
$\tau (U_e) \leq \tau(U).$
\item If $m$ is odd, then
$\tau(U_o) \leq \tau(U).$
\end{enumerate}
\end{thm}

Lov\'{a}sz\cite{Lovasz1972} proved that a graph of order $n$ with a unique perfect matching cannot have more than $\frac{n^2}{4}$ edges. An edge is said to be a \emph{cut edge} if its removal increases the number of components of a graph. Kotzig\cite{Kotzig1959} proved that a connected graph $G$ with a unique perfect matching contains a cut edge $uv$ that is an edge of the perfect matching of $G$. Based on the above two structural properties, Fan, Lin and Lu\cite{Fan-Lin-Lu} determined the graph attaining the maximum spectral radius among all graphs of order $2n$ with a unique perfect matching. Suppose that $G_1$ is an empty graph with vertex set $U=\{u_{1},u_{2},\ldots, u_{n}\}$, and $G_2$ is a complete graph with vertex set $W=\{w_{1},w_{2},\ldots, w_{n}\}$. Let $G(2n,1)$ be the graph of order $2n$ obtained from $G_{1}\cup G_{2}$ by letting $N_{G_2}(u_i)=\{w_{1},w_{2},\ldots,w_{i}\}$ for $1\leq i\leq n$. Clearly, $G(2n,1)$ contains a unique perfect matching.
\begin{thm}[Fan, Lin and Lu\cite{Fan-Lin-Lu}]\label{Fan-Lin-Lu-1}
If $G$ is a connected graph of order $2n$ with a unique perfect matching, then
$$\rho(G)\leq\rho(G(2n,1)),$$
with equality if and only if $G\cong G(2n,1)$.
\end{thm}

Observe that the Laplacian matrix of a disjoint union of $\frac{n}{2}$ edges has eigenvalues 0 and 2 in a graph $G$ of order $n$. This implies that deletion of the edges of a perfect matching in a graph $G$ will reduce the eigenvalues of the Laplacian matrix of $G$ with at most 2. Brouwer and Haemers\cite{A.B} proved that a regular graph of even order with algebraic connectivity $\mu_{n-1}$ has at least $\lfloor\frac{\mu_{n-1}+1}{2}\rfloor$ disjoint perfect matchings. Soon after, Cioab\u{a}, Gregory and Haemers\cite{S.C} considered the number of edge-disjoint matchings in a regular graph $G$ from the view of the second largest eigenvalue.

\begin{thm}[Cioab\u{a}, Gregory and Haemers\cite{S.C}]
 An $r$-regular graph $G$ of even order has at least $\lfloor\frac{r-\lambda_2(G)+1}{2}\rfloor$ edge-disjoint perfect matchings.
 \end{thm}
\subsection{Spanning trees}

In 1889, Cayley\cite{Cayley} showed that the number of spanning trees in the complete graph $K_n$ is $n^{n-2}$.
The edge-disjoint spanning trees has many applications in fault-tolerance networks as well as network reliability~\cite{Cunningham,Hobbs91}. Thus it is quite interesting to explore how many edge-disjoint spanning trees in a given graph. The \textit{spanning tree packing number} (or simply \textit{STP number}) of a graph $G$, denoted by $\tau(G)$, is the maximum number of edge-disjoint spanning trees contained in $G$. Nash-Williams~\cite{Nash-Williams} and Tutte~\cite{Tutte1961} independently discovered a fundamental theorem that characterizes graphs $G$ with $\tau(G)\ge k$.

For any partition $\pi$ of $V(G)$, let $E_{G}(\pi)$  denote the set of edges in $G$ whose endpoints lie in different parts of $\pi$, and let $e_{G}(\pi)=|E_{G}(\pi)|$.
\begin{thm}[Nash-Williams~\cite{Nash-Williams} and Tutte~\cite{Tutte1961}]
Let $G$ be a connected graph. Then $\tau(G)\geq k$ if and only if for any partition $\pi$ of $V(G)$, $$e_{G}(\pi)\geq k(t-1),$$
where $t$ is the number of parts in the partition $\pi$.
\end{thm}

The well-known Matrix-Tree Theorem of Kirchhoff \cite{Kirchhoff} indicates that the number of spanning trees (not necessarily edge-disjoint) of a graph $G$ with $n$ labelled vertices is $\frac{\prod_{i=1}^{n-1}\mu_{i}}{n}$. For edge-disjoint spanning trees, Seymour proposed the following problem in private communication to Cioab\u{a} as mentioned in~\cite{Cioaba-Wong}.
\begin{prob}\label{prob1}
Let $G$ be a connected graph. Determine the relationship between the spanning tree packing number $\tau(G)$ and the eigenvalues of $G$.
\end{prob}
Inspired by the Kirchhoff's Matrix-Tree Theorem and Problem~\ref{prob1}, Cioab\u{a} and Wong \cite{Cioaba-Wong} started to study the spanning tree packing number via the second largest eigenvalue of the adjacency matrix.
\begin{thm}[Cioab\u{a} and Wong~\cite{Cioaba-Wong}]
Let $k$ and $r$ be two integers with $r\geq 2k\geq 4$ and $G$ be a $r$-regular connected graph. If
 $$\lambda_{2}(G)<r-\frac{2(2k-1)}{r+1},$$
 then $\tau(G)\geq k$.
\end{thm}

Cioab\u{a} and Wong~\cite{Cioaba-Wong} further conjectured that the sufficient condition can be improved to $\lambda_{2}(G)<r-\frac{2k-1}{r+1}$. In the same paper, they did the preliminary work of this conjecture for $k=2,3$ and gave examples to show the bound is best possible. Later, Gu, Lai, Li and Yao\cite{Gu-Lai}, Li and Shi\cite{Li-Shi}, Liu, Hong and Lai\cite{Liu-Hong-Lai} independently extended the conjecture to graphs that are not necessarily regular.
\begin{conj}[See\cite{Gu-Lai,Li-Shi,Liu-Hong-Lai}]\label{Gu-Lai-Li-Yao-Hong}
Let $k\geq 2$ be an integer, and let $G$ be a graph with minimum degree $\delta\geq 2k$. If
$$\lambda_{2}(G)<\delta-\frac{2k-1}{\delta+1},$$
then $\tau(G)\geq k$.
\end{conj}

Gu, Lai, Li and Yao\cite{Gu-Lai} confirmed the Conjecture \ref{Gu-Lai-Li-Yao-Hong} for $k=2,3$ and obtained the following partial result for $k\geq 4$.
\begin{thm}[Gu, Lai, Li and Yao\cite{Gu-Lai}]
Let $k\geq 4$ be an integer, and let $G$ be a graph with minimum degree $\delta\geq 2k$. If
 $$\lambda_2(G)< \delta-\frac{3k-1}{\delta+1},$$
then $\tau(G)\geq k$.
\end{thm}


In \cite{Li-Shi}, Li and Shi obtained the following result, which suggests an approximate formulation of Conjecture \ref{Gu-Lai-Li-Yao-Hong} for a graph with a large order $n$ and small minimum degree $\delta$.
\begin{thm}[Li and Shi\cite{Li-Shi}]
Let $k \geq 2$ be an integer, and let $G$ be a graph of order $n$ with minimum degree $\delta \geq 2k$. If

 \begin{equation*}
\lambda_2(G)<\left\{
\begin{array}{ll}
\delta-\frac{n(k-1)}{(n-\delta-1)(\delta+1)} & \mbox{if $n \leq 3\delta+2$}, \\
\delta-\frac{8 k-7}{3(\delta+1)}  & \mbox{if $3(\delta+1) \leq n \leq 5(\delta+1)$}, \\
\delta-\frac{2 k-1}{\delta+1}+\frac{2(k-1)}{n-2(\delta+1)} & \mbox{if $n>5(\delta+1)$},
\end{array}
\right.
\end{equation*}
then $\tau(G)\geq k$.
\end{thm}

Liu, Hong and Lai\cite{Liu-Hong-Lai} made the investigation of Conjecture \ref{Gu-Lai-Li-Yao-Hong} and further proved that this conjecture is hold for sufficiently large $n$.
\begin{thm}[Liu, Hong and Lai\cite{Liu-Hong-Lai}]
Let $G$ be a graph of order $n \geq(2 k-1)(\delta+1)$ with minimum degree $\delta \geq 2k \geq 4$. If
$$
\lambda_2(G) \leq \delta-\frac{2 k-2 / k}{\delta+1} ~\text { or }~\lambda_2(G) \leq \delta-\frac{2 k-1}{\delta+1}\text {,}
$$
 then $\tau(G) \geq k$.
\end{thm}
Conjecture \ref{Gu-Lai-Li-Yao-Hong} was completely settled in 2014 by Liu et al.~\cite{Liu-Hong} who proved a stronger result, which also implied the truth of the conjecture of Cioab\u{a} and Wong \cite{Cioaba-Wong}. Moreover, the bound in Conjecture \ref{Gu-Lai-Li-Yao-Hong} has been shown to be essentially best possible in \cite{coppww22} by constructing extremal graphs.

For given integers $\delta$ and $g$ with $\delta>0$ and $g \geq 3$, let $t=\lfloor\frac{g-1}{2}\rfloor$, and define the Moore function as follows.
$$
f(\delta, g)= \begin{cases}2 t+1 & \text { if } \delta=2 \text { and } g=2 t+1 \\ 1+\delta+\sum_{i=2}^t(\delta-1)^i & \text { if } \delta \geq 3 \text { and } g=2 t+1 \\ 2 t+2 & \text { if } \delta=2 \text { and } g=2 t+2 \\ 2+2(\delta-1)^t+\sum_{i=1}^{t-1}(\delta-1)^i & \text { if } \delta \geq 3 \text { and } g=2 t+2\end{cases}
$$
By utilizing the Moore function, Liu, Lai, and Tian\cite{Liu-Lai-Tian} investigated the spanning tree packing number of a graph. If a graph $G$ has a cycle, the girth of $G$ is
the length of the shortest cycle in $G$.
\begin{thm}[Liu, Lai and Tian\cite{Liu-Lai-Tian}]
Let $g$ and $k$ be integers with $g \geq 3$ and $k \geq 2$, and $G$ be a graph of order $n$ with minimum degree $\delta \geq 2 k \geq 4$ and girth $g$.
 If
$$\lambda_2(G)<\delta-\frac{2 k-1}{f(\delta, g)},$$
  then $\tau(G) \geq k$.
\end{thm}

Let $\mathcal{G}_t$ be the set of graphs such that for each graph $G \in \mathcal{G}_t$ there exist at least $t+1$ non-empty disjoint proper subsets $V_1, V_2, \ldots, V_{t+1} \subseteq V(G)$ satisfying $V(G) \backslash(V_1 \cup V_2 \cup \cdots \cup V_{t+1}) \neq \varnothing$ and edge connectivity $\kappa^{\prime}(G)=e_{G}(V_i, V(G) \backslash V_i)$ for $i=1,2, \ldots, t+1$.
 Duan, Wang and Liu\cite{Duan-Wang-Liu} gave the following spectral condition to guarantee $\tau(G) \geq k$ in terms of its third largest eigenvalue among all $\mathcal{G}_1$.
\begin{thm}[Duan, Wang and Liu\cite{Duan-Wang-Liu}]\label{Duan-Wang-Liu}
Let $k \geq 2$ be an integer, and let $G \in \mathcal{G}_1$ be a graph with minimum degree $\delta \geq 2 k$ and maximum degree $\Delta$. If
 $$\lambda_3(G)<2 \delta-\Delta-\frac{2(3 k-1) }{\delta+1},$$
 then $\tau(G) \geq k$.
\end{thm}

Later, Hu, Wang and Duan\cite{Hu-Wang-Duan} improved the upper bound in Theorem \ref{Duan-Wang-Liu} to $\delta-\frac{8(2k-1)}{3(\delta+1)}$, and determined the relationship between $\lambda_4(G)$ and spanning tree packing number of a graph $G \in \mathcal{G}_2$.

\begin{thm}[Hu, Wang and Duan\cite{Hu-Wang-Duan}]
Let $k \geq 2$ be an integer, and let $G \in \mathcal{G}_2$ be a graph with minimum degree $\delta \geq 3 k$. If $$\lambda_4(G)<\delta-\frac{9 k-3}{\delta+1},$$
then $\tau(G) \geq k$.
\end{thm}

 Motivated by the above results, Fan, Gu and Lin\cite{Fan-Gu-Lin} studied the spanning tree packing number by means of the spectral radius of graphs. They first investigated an edge extremal result for $\tau(G)\ge k$.

\begin{thm}[Fan, Gu and Lin\cite{Fan-Gu-Lin}]\label{thm::edgenumber}
Let $G$ be a connected graph with minimum degree $\delta\geq 2k$ and order $n\geq 2\delta+2$. If
$$e(G)\geq {\delta+1\choose2}+{n-\delta-1\choose 2} +k,$$
then $\tau(G)\geq k$.
\end{thm}
The condition in Theorem~\ref{thm::edgenumber} is tight. Denote by $K_n$ the complete graph on $n$ vertices, and $\mathcal{G}_{n,n_1}^{i}$ the set of graphs obtained from $K_{n_1}\cup K_{n-n_{1}}$ by adding $i$ edges between $K_{n_1}$ and $K_{n-n_{1}}$. Notice that any graph $G$ in $\mathcal{G}_{n,\delta+1}^{k-1}$ has exactly ${\delta+1\choose2} + {n-\delta-1\choose 2} +k-1$ edges but $\tau(G)<k$. In the same paper, they then focus on a spectral analogue. Let $B_{n,\delta+1}^{i}$ be a graph obtained from $K_{\delta+1}\cup K_{n-\delta-1}$ by adding $i$ edges joining a vertex in $K_{\delta+1}$ and $i$ vertices in $K_{n-\delta-1}$. They discovered a sufficient condition for $\tau(G)\geq k$ via the spectral radius, and characterized the unique spectral extremal graph $B_{n,\delta+1}^{k-1}$ among the structural extremal graph family $\mathcal{G}_{n,\delta+1}^{k-1}$.
\begin{thm}[Fan, Gu and Lin\cite{Fan-Gu-Lin}]\label{Fan-Gu-Lin}
Let $k\geq 2$ be an integer, and let $G$ be a connected graph with minimum degree $\delta\geq 2k$ and order $n\geq 2\delta+3$. If $\rho(G)\geq \rho(B_{n,\delta+1}^{k-1})$, then $\tau(G)\geq k$ unless $G\cong B_{n,\delta+1}^{k-1}$.
\end{thm}

Let
$$\eta(G)=\min\Big\{\frac{|X|}{\omega(G-X)-\omega(G)}\Big\}$$
where the minimum is taken over all edge subsets $X$ in $G$. Nash-Williams-Tutte Theorem indicates that for a connected graph $G$, $\tau(G) \geq k$ if and only if $\eta(G) \geq k$. Since $\eta(G)$ is possibly fractional, we have $\tau(G)=\lfloor\eta(G)\rfloor$. Thus, we call $\eta(G)$ is the fractional spanning tree packing number of $G$. Hong, Gu, Lai and Liu\cite{Hong-Gu-Lai-Liu} investigated the relationship between $\eta(G)$ and the eigenvalue of $G$.
\begin{thm}[Hong, Gu, Lai and Liu\cite{Hong-Gu-Lai-Liu}]\label{Hong-Gu-Lai-Liu}
 Let $s,t$ be two positive integers, and let $G$ be a graph with minimum degree $\delta \geq \frac{2 s}{t}$. If
 $$\lambda_2(G)<\delta-\frac{2 s-1}{t(\delta+1)},$$
 then $\eta(G) \geq \frac{s}{t}$.
\end{thm}
When $s/t$ is a positive integer, this result also confirms the Conjecture \ref{Gu-Lai-Li-Yao-Hong}.
\section{Some possible problems}

This is a temporary section for discussing and some listing related problems, with a particular focus on spectral radius of graphs.

The existence of rainbow matchings in graphs has been extensively investigated by many researchers in past several decades.
 Joos and Kim \cite{JK20} posted the following question:
 \begin{prob}[Joos and Kim \cite{JK20}]\label{Joos-Kim}
 Let $H$ be a graph with $t$ edges, $\mathbf{G}$ be a family of graphs and $\mathcal{G}=\{G_1,\ldots,G_t\}$ be a collection of not necessarily distinct graphs on the same vertex set $V$ such that $G_i\in \mathbf{G}$ for all $1\leq i\leq t$. Which properties imposed on $\mathbf{G}$ yield a $\mathcal{G}$-rainbow graph isomorphic to $H$?
\end{prob}
In the same paper, they answered the Problem \ref{Joos-Kim} from version of Dirac's theorem.
\begin{thm}[Joos and Kim \cite{JK20}]
 Let $\mathcal{G}=\{G_1,\ldots,G_{n/2}\}$ be a collection of  graphs on the same vertex set $[n]$ such that $\delta(G_i)\geq n/2$ for $1\leq i\leq n/2$. Then $\mathcal{G}$ admits a rainbow matching.
\end{thm}
Specially, we would like to propose the following problems from the perspective of edge and spectral radius, respectively.
\begin{prob}
Let $n,r$ be two integers such that $1\leq r< n/2$. Let $\mathcal{G}=\{G_1,\ldots,G_{n/2}\}$ be a collection of  graphs on the same vertex set $[n]$ such that $\delta(G_i)\geq r$ and $e(G_i)>e(K_{r-1}\vee (rK_1\cup K_{n-2r+1}))$ for $1\leq i\leq n/2$, then $\mathcal{G}$ admits a rainbow matching.
\end{prob}
\begin{prob}
Let $n,r$ be two integers such that $1\leq r< n/2$. Let $\mathcal{G}=\{G_1,\ldots,G_{n/2}\}$ be a collection of  graphs on the same vertex set $[n]$ such that $\delta(G_i)\geq r$ and $\rho(G_i)>\rho(K_{r-1}\vee (rK_1\cup K_{n-2r+1}))$ for $1\leq i\leq n/2$, then $\mathcal{G}$ admits a rainbow matching.
\end{prob}

For the existence of rainbow matchings in bipartite graphs, Aharoni and Howard \cite{AH11,AH16} posed a conjecture.

\begin{conj}[Aharoni and Howard \cite{AH11,AH16}]\label{AH11,AH16}
 Let $\mathcal{G}=\{G_1,\ldots,G_{t}\}$ be $t$ bipartite graphs such that $\Delta(G_i)\leq \Delta$ and $e(G_i)>(t-1)\Delta$. Then $\mathcal{G}$ admits a rainbow matching.
\end{conj}

By incorporating rainbow matching and spectral radius conditions in a bipartite graph, we can consider the above Conjecture \ref{AH11,AH16} from another perspective.
\begin{prob}
 Let $\mathcal{G}=\{G_1,\ldots,G_{t}\}$ be $t$ bipartite graphs such that $\Delta(G_i)\leq \Delta$. Find a spectral radius condition for a collection of graphs $\mathcal{G}$ to admits a rainbow matching.
\end{prob}

An edge-colored graph is called \textit{rainbow} if the colors on its edges are distinct. For graphs $G$ and $H$, the \textit{anti-Ramsey number} is the maximum number of colors in an edge-colored $G$ with no rainbow copy of $H$. Note that the rainbow number is closely related to anti-Ramsey numbers.
Then we can considered the following problem.

\begin{prob}
Find a spectral condition for rainbow matching in edge-colored graphs (anti-Ramsey number spectral version).
\end{prob}

Concerning connected factors in a graph, some researchers have focused on finding spectral conditions for a graph to contain a Hamiltonian cycle or Hamiltonian path. In \cite{O2022}, O raised a question regarding the identification of specific eigenvalue conditions that guarantee the existence of a connected (even or odd) $[a,b]$-factor in an $h$-edge-connected $r$-regular graph. It is therefore natural to inquire about the corresponding question for a more general condition.

\begin{prob}
Characterize the eigenvalue conditions for the existence of a connected $[a,b]$-factor in graphs.
\end{prob}



For two positive integers $n$ and $k$, let $m(2n,k)$ be the maximum number of edges in a graph of order $2n$ with a unique $k$-factor. In 1984, Hendry\cite{Hendry} proved $m(2n,2)=\lfloor\frac{n(2n+1)}{2}\rfloor$ and conjectured $m(2n,k)=nk+\binom{2n-k}{2}$ for $k> n$, and $kn$ is even and $m(2n,k)=n^2+\frac{(k-1)n}{2}$ for $n=ks$, where $s$ is a positive integer. In 2000, Johann\cite{Johann} confirmed this conjecture and obtained the corresponding extremal graphs $K_{2n-k}\nabla H$, where $H$ is a $2(k-n)$-regular graph, for $k> n$ and $G(2n,k)$ for $n=ks$, respectively. The construction of $G(2n,k)$ is shown below.

Let $n=sk+t$ with $s\geq 1$ and $0\leq t\leq k-1$. We give the process to construct the graph $G(2n,k)$. First define a graph $F_1$ on $2(k+t)$ vertices as follows. Let $H_{1}\cong K_{t}\nabla tK_1$, and let $A_{11}=V(tK_{1})$ and $A_{12}=V(K_t)$. Denote by $H_2$ the graph obtained from $kK_{1}\cup K_{k}$ by adding edges between $V(kK_{1})$ and $V(K_{k})$ such that $d_{H_2}(v)=k-t$ for $v\in V(kK_{1})$ and $d_{H_2}(u)=2k-t-1$ for $u\in V(K_{k})$. Suppose that $A_{21}=V(kK_{1})$ and $A_{22}=V(K_{k})$. Let $F_1$ be the graph of order $2(k+t)$ obtained from $H_{1}\cup H_{2}$ by connecting all vertices of $A_{1j}$ with $A_{2(3-j)}$ for $1\leq j\leq 2$ and adding all edges between $A_{12}$ and $A_{22}$. The resulting graph $F_1$ has exactly one $k$-factor. Suppose that $U_1=A_{11}\cup A_{21}$ and $W_1=A_{12}\cup A_{22}$. Next take $s-1$ copies of $K_{k}\nabla kK_1$ labeled $F_{2},\ldots, F_{s}$. For $2\leq i\leq s$, let $U_{i}$ and $W_{i}$ be the vertices set of $V(kK_1)$ and $V(K_k)$ in each $F_i$, respectively. Then the graph $G(2n,k)$ is obtained by adding edges connecting all vertices of $W_i$ in $F_i$ to all vertices in $F_j$ for each $i,j$ with $1\leq i< j\leq s$. The resulting graph $G(2n,k)$ has a unique $k$-factor.

From the perspective of spectral radius condition, we have the following problem.

\begin{prob}\label{prob::5.1}
What is the maximum spectral radius and what is the corresponding extremal graph among all graphs with a unique $k$-factor for $k\geq 1$?
\end{prob}

In this paper, Theorem \ref{Fan-Lin-Lu-1} has given the answer to Problem \ref{prob::5.1} for $k=1$. It is $G(2n,1)$. However, the structure of graphs with a unique $k$-factor is more complicated for $k\geq 2$, and it seems difficult to determine the extremal graphs about the problem.

\begin{prob}
For $k\geq 2$. Suppose that $G$ is a graph of order $2n$ with a unique $k$-factor.

\noindent (I) Does $\rho(G)\leq\rho(K_{2n-k}\nabla H)$, where $H$ is a $2(k-n)$-regular graph for $k> n$?

\noindent (II) Does $\rho(G)\leq\rho (G(2n,k))$ for $k\leq n$?
\end{prob}

Theorem~\ref{Fan-Gu-Lin} actually implies that $B_{n,\delta+1}^{\tau}$ is the unique graph that has the maximum spectral radius among all graphs of fixed order $n$ with minimum degree $\delta$ and spanning tree packing number $\tau$. Let $G$ be a minimal graph with $\tau(G)\ge k$ and of order $n$, that is, $G$ consists of exactly $k$ edge-disjoint spanning trees with no extra edges. This implies that $\tau(G)=k$ and $e(G)=k(n-1)$. We are interested in the maximum possible spectral radius of $G$.

\begin{prob}\label{prob::minkst}
Let $G$ be a minimum graph with $\tau(G)\ge k$ and of order $n$. For each $n\ge 4$ and each $k\ge 2$, determine the maximum possible spectral radius of $G$ and characterize extremal graphs.
\end{prob}

The following sharp upper bound on the spectral radius was obtained by Hong, Shu and Fang\cite{Hong-Shu-Fang} and Nikiforov\cite{Nikiforov2002}, independently.
\begin{thm}[See \cite{Hong-Shu-Fang,Nikiforov2002}]\label{HSF-N}
Let $G$ be a graph on $n$ vertices and $m$ edges with minimum degree $\delta\geq 1$. Then
$$\rho(G) \leq \frac{\delta-1}{2}+\sqrt{2 m-n \delta+\frac{(\delta+1)^{2}}{4}},$$
with equality if and only if $G$ is either a $\delta$-regular graph or a bidegreed graph
in which each vertex is of degree either $\delta$ or $n-1$.
\end{thm}

To attack Problem~\ref{prob::minkst}, notice that $\delta(G)\ge k$ and $e(G)=k(n-1)$, we can easily obtain an upper bound on $\rho(G)$ by Theorem~\ref{HSF-N}. However, this upper bound may not be tight, since for most values of $n$, $G$ cannot be $k$-regular or a bidegreed graph in which each vertex is of degree either $k$ or $n-1$. To attain the maximum spectral radius, it seems that $G$ is obtained from a $K_{2k}$ by continuously adding a vertex and $k$ incident edges step by step until $G$ has $n$ vertices, and in particular, we guess the graph is $K_{k}\nabla (K_{k}\cup (n-2k)K_1)$.

As a dual problem of spanning tree packing, Nash-Williams~\cite{Nash64} ever studied the forest covering problem, seeking the minimum number of forests that cover the entire graph.  The {\it arboricity} $a(G)$ is the minimum number of edge-disjoint forests whose union equals $E(G)$.
\begin{thm}[Nash-Williams~\cite{Nash64}]\label{thm::5.1}
Let $G$ be a connected graph. Then $a(G)\le k$ if and only if for any subgraph $H$ of $G$, $|E(H)|\le k(|V(H)|-1)$.
\end{thm}

Naturally, we have the following problem.
\begin{prob}\label{prob::arbo}
Find a tight spectral radius condition for a graph $G$ of order $n$ with $a(G)\le k$ and characterize extremal graphs.
\end{prob}

When $n\leq 2k$, Problem~\ref{prob::arbo} is trivial. Any graph $G$ of order $n\leq 2k$ has the property $a(G)\le k$ and so the extremal graph is $K_n$. To see this, for any subgraph $H$ of $G$, we can deduce that $|E(H)|\leq {|V(H)|\choose 2} \leq k(|V(H)|-1)$ when $n\leq 2k$, and the conclusion follows from Theorem~\ref{thm::5.1}. The situation becomes more involved for $n\geq 2k+1$. In fact, this case is even stronger than Problem~\ref{prob::minkst}, since if $G$ consists of exactly $k$ edge-disjoint spanning trees with no extra edges, we have $a(G) =\tau(G)=k$.
Notice that $a(K_{k}\nabla (K_{k}\cup (n-2k)K_1))=k$ and $e(K_{k}\nabla (K_{k}\cup (n-2k)K_1))=k(n-1)$, thus we guess that the extremal graph w.r.t. the spectral radius is also $K_{k}\nabla (K_{k}\cup (n-2k)K_1)$.

The eigenvalues of a graph are closely interconnected with graph parameters, such as chromatic number\cite{Nikiforov2002}, clique number\cite{Wilf1986} and cut edge\cite{Liu-Lu-Tian}. For more results, we refer the reader to \cite{Haemers,Hoffman}. It is worth noting that the majority of researchers have focused primarily on studying the maximum spectral radius of a graph under fixed conditions of a single graph parameter, with relatively less attention paid to the results obtained from giving two graph parameters.
Therefore, it is interesting to consider the maximum spectral radius and corresponding extremal graphs among all graphs with two given parameters.
\begin{prob}
Find a tight spectral radius condition for a graph $G$ of order $n$ with $\Delta(G)\leq r$ and $\beta(G)\leq s$.
\end{prob}

\end{document}